\documentclass[dvips,11pt]{amsart}
\usepackage{graphicx}
\usepackage{psfrag}
\usepackage{amsfonts}
\usepackage{amsmath}
\usepackage{amsthm}
\usepackage{amssymb}
\usepackage{stmaryrd}
\usepackage{bbold}
\theoremstyle{plain}
\newtheorem{thrm}{Theorem}[section]
\newtheorem{lmm}[thrm]{Lemma}

\newtheorem{lgrthm}[thrm]{Algorithm}
\theoremstyle{definition}

\newtheorem{rmrk}[thrm]{Remark}

\newcommand{\virg}{\;,\;\;}
\newcommand{\ee}{\mathbf{e}}
\newcommand{\xx}{\mathbf{x}}
\newcommand{\yy}{\mathbf{y}}
\newcommand{\VV}{\mathbf{V}}
\newcommand{\FF}{\mathbf{F}}
\newcommand{\ff}{\mathbf{f}}
\newcommand{\uu}{\mathbf{u}}
\newcommand{\p}{\text{p}}
\newcommand{\IR}{\mathbb{R}}
\newcommand{\wsto}{\stackrel{\star}{\relbar\joinrel\rightharpoonup}}
\newcommand{\GGij}{\mathbf{G}_{ij}}
\newcommand{\GGik}{\mathbf{G}_{ik}}
\newcommand{\ggamma}{{\boldsymbol\gamma}}
\newcommand{\gammaij}{{\gamma}_{ij}}
\newcommand{\llambda}{{\boldsymbol\lambda}}
\newcommand{\lambdaij}{{\lambda}_{ij}}

\newcommand{\vseq}{\vspace{6pt}}

\newcommand{\findemo}{\,\,\,$\Box$\vspace{0.5cm}}

\begin{document}
\title{Numerical Simulation of Gluey Particles}
\author{Aline Lefebvre$^1$}
\thanks{1. Laboratoire de Math\'ematiques, Universit\'e Paris-Sud, 91405 Orsay Cedex, France, aline.lefebvre@math.u-psud.fr}


\keywords{fluid/particle systems, fluid/solid interaction,
  lubrication force, contacts, Stokes fluid}

\begin{abstract}
We propose here a model and a numerical scheme to compute the motion
of rigid particles interacting through the lubrication force. In the
case of a particle approaching a plane, we propose an algorithm and
prove its convergence towards the solutions to the gluey particle model
described in~\cite{Maury2007}. We propose a multi-particle version of
this gluey model which is based on the projection of the velocities
onto a set of admissible velocities. Then, we describe a multi-particle algorithm
for the simulation of such systems and present numerical results.

\bigskip
\noindent{\sc R\'esum\'e. } 
Nous proposons ici un mod\`ele ainsi qu'un sch\'ema num\'erique afin
de r\'esoudre le mouvement de particules rigides en interaction \`a
travers la force de lubrification. Dans le cas d'une particule \`a
l'approche d'un plan, nous proposons un algorithme et montrons sa
convergence vers le mod\`ele de particules visqueuses d\'ecrit
dans~\cite{Maury2007}. Nous proposons une version multi-particules de
ce mod\`ele qui est bas\'ee sur la projection des vitesses sur un
espace de vitesses admissibles. Ensuite, nous d\'ecrivons un
algorithme multi-particules pour la simulation de tels syst\`emes et pr\'esentons des
r\'esultats num\'eriques.
\end{abstract}

 \maketitle

\section{Introduction}

Slurries, lava's flows or red cells in blood are systems made of rigid
particles embedded in viscous fluids (if we consider as a first
approximation that red blood cells are rigid). Such systems can also
be found in industry: concrete, paper pulp
or some food industry products. These systems present varieties of
noticeable rheological behaviours, whose study has been the subject of a great amount
of researches, with contributions coming from engineering, chemistry, physics
or mathematics. The basic problem is to predict macroscopic transport properties
of these suspensions -- viscosity, settling velocity -- from
microstructures, that is to say, from the interactions between
particles and from their spatial configuration.\vseq

In case of dilute suspensions, theoretical results come from
neglecting near field interactions. For example, in 1906, Einstein
proposed an asymptotic formula for the apparent viscosity of dilute
suspensions~\cite{Einstein1906}. In that case, apparent viscosity only
depends on the solid volume fraction. Unfortunately, agreement
between such asymptotic results and experiments generally fails as soon
as the solid fraction reaches a few percent. For higher solid
fractions, near field interactions can not be neglected anymore and it
becomes essential to take them into account. 
Note that, studying the
behaviour of neighbouring particles is of great interest, not only to 
understand the behaviour of dense suspensions, but also to
study the fluid/particle system of equations modelling suspensions of
particles. Indeed, existence of solutions to these equations has been
proved as long as the distance remains stricly positive (see for
example~\cite{Desjardins1999,Takahashi2003,Takahashi2003b}). Global
weak solutions have also been constructed in~\cite{Feireisl2003,SanMartin2002}, supposing that
solids stick after contact. However, nothing is said concerning the
possibility that such a contact may occur in finite time. A good understanding of near field
interactions is therefore necessary to study more precisely these systems. \vseq

These interactions between solids embedded in a viscous fluid are due
to lubrication forces: for the solids to get very close, the fluid
must be evacuated from the narrow gap between them, which creates a force penalizing their relative
motion. This force is singular in the distance and this singularity is
sufficient to avoid contacts. Indeed, it has been proved
in~\cite{Hill} that in two dimensions, a smooth particle embedded in a
viscous fluid following Navier-Stokes equations can not touch a plane
in finite time. This behaviour can be recovered from the asymptotic
expansion of the lubrication force, available for a Stokes fluid
in three dimensions (see~\cite{Cox1967} for example):
\begin{equation}
\label{DL_Flub_1}
\FF_{lub}\sim -6\pi\mu r^2 \frac{\dot q}{q},
\end{equation}
where $\mu$ is the viscosity of the fluid, $r$ the radius of the
particle and $q$ the distance between the particle and the
plane. Indeed, using this first order approximation, we can write the Fundamental Principle of
Dynamics for a particle of mass $m$ submitted to an
external force $f$:
\begin{equation}
\label{PFD}
m\ddot q(t) = -6\pi\mu r^2 \frac{\dot q}{q} + mf(t),
\end{equation}
and the fact that the maximal solution to this ODE is global and never goes to
zero (contact) in finite time comes from the Cauchy-Lipschitz
theorem. Similarly, in case of a fixed sphere of radius $r_1$ and
another sphere of radius $r_2$ moving at velocity $\VV$ along the axe
of the centers, the first term of the developpement of the lubrication
force exerted on the moving particle is (see~\cite{Cox1974}):
\begin{equation}
\label{DL_Flub_2}
\FF_{lub}\sim -6\pi\mu
\frac{r_1^2r_2^2}{(r_1+r_2)^2}\frac{\VV}{q},
\end{equation}
and no contact can occur in finite time.\vseq

This force, while acting at microscopic level, can be very
important for the macroscopic behaviour of the global system,
especially in case of high density of particles. Even for
Stokes flows, it induces complexity
and nonlinearity. This complex link between microscopic and
macrosopic levels makes it
difficult to obtain theoretical results and studying these systems
requires numerical simulations. In order to obtain relevant
simulations for dense suspensions, the lubrication force has to be
taken into account with accuracy. However, direct numerical simulation induces space
discretization which makes it difficult to solve accurately the fluid
in the narrow gap between neighbouring particles. As a consequence,
numerical contacts can be observed in such simulations and physical
reasons as well as numerical robustness make it necessary to
develop specific technics to deal with these contacts.\vseq

A first idea to solve this problem is to search for a strategy allowing
an accurate computation of the lubrication forces. In~\cite{Hu1996}, a
method based on local refinements of the space and time meshes is
proposed, so that the lubrication force in the interparticle gap is
taken into account with accuracy and prevents overlappings. However, the number of refinements
needed is not known a priori and the method can become computationally
heavy. Consequently, less time-consuming methods have been
developped. Some of them consist in adding a short range repulsive
force (see~\cite{Glo1999,Glo2000} or~\cite{Turek2006}). In~\cite{Maury1999} a
minimizing algorithm is used to impose a minimal distance between the
particles, while in~\cite{Singh2003}, the particles are allowed to undergo
slight overlappings and an elastic repulsive force is added when such
overlappings are detected. All of these methods ensure numerical
robustness but introduce new parameters and do not take into account
the underlying physics. Another approach is to use inelastic
collisions. This idea has been proposed in~\cite{Tez1996} in order to impose a
minimal distance between the particles. In~\cite{Maury2006}, a scheme for
inelastic collisions, based on a global projection step of the
velocities, has been developped for granular flows and makes it
possible to handle lots of particles. This scheme has been coupled
with a fluid/particle solver in~\cite{GTN}, to avoid contacts. More physical strategies, taking
the lubrication force into account, have finally been proposed. Each
of them relies on the asymptotic developpement of the lubrication
force~(\ref{DL_Flub_2}). In~\cite{Dance2003,Nasseri2000}, it is shown that these lubrication
forces are solution to a linear system. They are
computed at each time step and added to the
simulations. Unfortunately, this leads to stiff systems and, whereas
it better takes into account the underlying physics, contact problems still occur because
of the time discretization. In~\cite{Maury1997} a method is proposed to
stabilize this problem by computing accurately sensible quantities
such as the interparticle distances. However, a projection step is
still needed for big time steps, in order to avoid overlappings.\vseq

The purpose of this article is to propose a strategy dealing simultaneously with contacts and
lubrication forces. We restrict ourself here to the study of a gluey
contact model without taking the surrounding fluid into account. This
model is based on  the gluey particle one
described in~\cite{Maury2007}. We propose an algorithm for this
particle/plane model and prove its convergence. Then, we generalize
it to the multi-particle case. The numerical strategy is to
combine the algorithm given for the plane/particle case with the
scheme proposed in~\cite{Maury2006} for granular flows. While programming
this multi-particle algorithm, we watched out for dealing with contacts
efficiently in order to manage to simulate collections of many
particles. Numerical simulations for few thousands of gluey particles are
presented in the last section. An example of coupling with a fluid/particle solver is given in
section~\ref{sec:couplage} in the particle/plane case. 

\section{Single particle above a plane}
\label{sec:plane_part}

\subsection{The gluey particle model}

\begin{figure}[hbtp]
\psfragscanon
\psfrag{x}[l]{$\ee_1$}
\psfrag{y}[l]{$\ee_2$}
\psfrag{U}[l]{$\VV=u\ee_2$}
\psfrag{q}[l]{$q$}
\psfrag{a}[l]{$r$}
\centering
\includegraphics[width=0.3\textwidth]{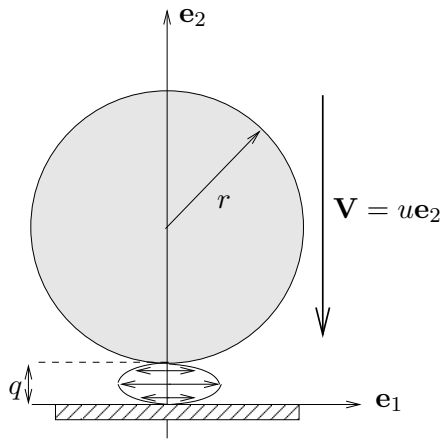}
\caption{Notations.}
\label{Notations}
\end{figure}

We consider a three-dimensional spherical particle moving perpendicularly to a plane (See Fig.~\ref{Notations}). Its
velocity and radius are denoted by $\VV$ and $r$ respectively. Its
distance to the plane is $q$.

The gluey particle model has been proposed in~\cite{Maury2007}. It
describes, from a macroscopic point of view, the behaviour of the
system near contact. It is built as the vanishing viscosity limit of the
lubrication model~(\ref{DL_Flub_1}) and relies on two states, glued ($q=0$) or unglued
($q>0$). These states are described by a new variable $\gamma$ which
stands for an adhesion potential: the more $\gamma$ is negative, the
more the solids are glued. 

We denote by $I=]0,T[$ the time interval. The unknowns $q$ and $\gamma$ belong to the
following functional spaces:
\begin{equation*}
q\in W^{1,\infty}(I)\virg \dot q\in BV(I)\virg\gamma\in
BV(I)\hbox{,}
\end{equation*}
and the initial conditions are:
\begin{equation*}
q(0)=q^0>0\virg \dot q (0) =u^0\virg \gamma(0)=0.\vseq
\end{equation*}
In order to be able to generalize the model to the multi-particle
case, we use the following second order ODE formulation given in~\cite{Maury2007}:
\begin{equation}
\dot q(t^+) = \Pi_{C_{q,\gamma}(t)}\dot q(t^-)\label{loi_choc},\vseq
\end{equation}
\begin{equation}
m \ddot q =  mf + \lambda \hbox{ in } {\mathcal M}(I)=\left({\mathcal C}_c (I) \right)'\label{PFD_contraint},\vseq
\end{equation}
\begin{equation}
\hbox{supp}(\lambda)\subset \{t\virg q(t)=0\}\label{supp_lambda},\vseq
\end{equation}
\begin{equation}
\dot \gamma = -\lambda\label{evol_gamma},\vseq
\end{equation}
\begin{equation}
q\geq 0 \virg \gamma\leq 0\label{contrainte_q_gamma},\vseq
\end{equation}
where $C_{q,\gamma}(t)$ is the set of admissible velocities at time $t$:
$$
C_{q,\gamma}(t)=\left|
\begin{array}{ll}
\{0\} & \hbox{ if }  \, \gamma(t^-)<0, \vseq\\
\IR^+ &  \hbox{ if } \, \gamma(t^-)=0  \virg q(t)=0,\vseq \\
\IR & \hbox{ else }.
\end{array}
\right.
$$
\begin{rmrk}
In this formulation, $\dot q$ and $\gamma$ are supposed to be in $BV(I)$. In order to
  alleviate the notations, their differential measures have been
  denoted by $\ddot q$ and $\dot \gamma$ respectively.
\end{rmrk}

\begin{figure}[hbtp]
\psfrag{q(t)}[l]{\scriptsize $ q(t)$}
\psfrag{gamma(t)}[l]{\scriptsize  $ \gamma (t)$}
\psfrag{f=-2}[l]{\scriptsize  $ f=-2$}
\psfrag{f=2}[l]{\scriptsize  $ f=2$}
\psfrag{ddqf}{\scriptsize  $ \ddot q = f$}
\psfrag{dgf}{\scriptsize  $ \dot\gamma = mf$}
\psfrag{q0}{\scriptsize  $ q=0$}
\psfrag{timp}{\scriptsize  $ t_{hitting}$}
\psfrag{gtimp}{\scriptsize  $ \gamma(t_{hitting}^+)=m\dot q(t_{hitting}^-)$}
\centering
\includegraphics[width=0.7\textwidth]{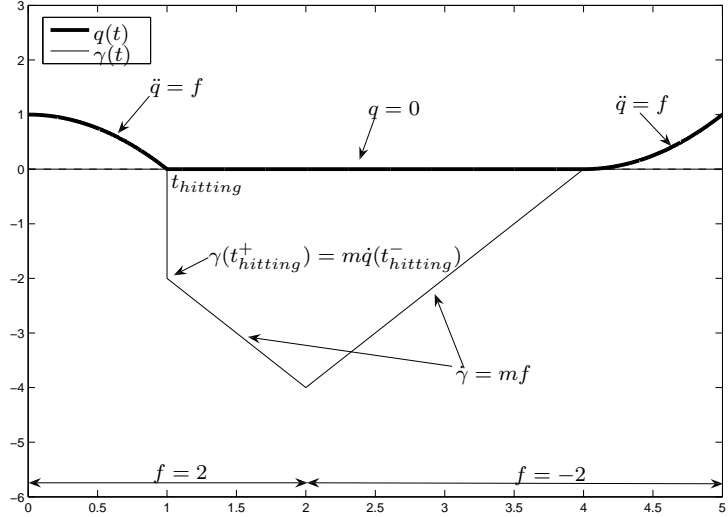}
\caption{Example of solution to the gluey particle model.}
\label{schema_modele_limite}
\end{figure}

The behaviour of the solutions to this problem is the
following. By~(\ref{PFD_contraint}) and~(\ref{supp_lambda}), $q$ is
solution to $\ddot q=f$ while there is no contact ($q>0$).  Suppose a collision occurs
at time $t_0$, we have $\dot q (t^-_0)<0$ and $\gamma(t^-_0)=0$. Then
$C_{q,\gamma}(t_0)$ is $\IR^+$ and~(\ref{loi_choc}) gives $\dot
q(t^+_0)=0$. By~(\ref{PFD_contraint}), we obtain that, in the sense of
distributions, $\lambda$ identifies to the dirac mass at time $t_0$
weighted by the velocity jump $m(\dot q(t^+_0)-\dot
q(t^-_0))=-m\dot q(t^-_0)$. This, together with~(\ref{evol_gamma}) finally
gives that $\gamma$ is initialized to the
value $m\dot q(t^-_0)<0$. From then, while $\gamma$ remains strictly
negative, $C_{q,\gamma}$ is reduced to $\{0\}$ and, combining this
with~(\ref{loi_choc}) gives that there is adhesion between the
solids ($q=0$). During this adhesion, $\ddot q$ is zero and therefore,~(\ref{evol_gamma}) associated
to~(\ref{PFD_contraint}) gives $\dot\gamma=mf$. By definition of $C_{q,\gamma}$, the particle is
allowed to take off when $\gamma$ is back to zero. An example of such a behaviour is given
in figure~\ref{schema_modele_limite}

\begin{rmrk}
The additional constraints~(\ref{contrainte_q_gamma}) are
necessary. Indeed, suppose that $t_1$ is
an unsticky contact time ($q(t_1)=0$, $\gamma(t_1^-)=0$, $\dot q
(t_1^-)=0$). If the force is negative after this instant and if
we do not impose $q\geq 0$, then $\gamma\equiv 0$ and $\dot
q(t)=\int_{t_1}^t f(s)ds$ is a solution to the problem and the particle can enter the
wall. Similarly, if the force is positive and if the constraint
$\gamma\geq 0$ is not imposed, $q\equiv 0$ and $\gamma(t)=m\int_{t_1}^t
f(s)ds$ is a solution and $\gamma$ can become strictly positive.
\end{rmrk}

Before proposing an algorithm to compute the solutions to this model,
we make a few remarks about its interpretation and its physical relevance.

\begin{rmrk}[Physical interpretation]
As already mentioned, a smooth particle embedded in a newtonian fluid never touches the plane in
finite time. In the context of the gluey particle model, the variable $q$ can be seen as a macroscopic
distance between the solids: it is equal to zero as soon as the solids
are near contact. The new variable $\gamma$, which is obtained as the
limit of $\gamma_\mu=\mu ln(q)$, stands
for the microscopic distance.
To understand the behaviour of the gluey particle system presented on
figure~\ref{schema_modele_limite}, one can consider a rigid ball falling on a table
coatted with a viscous fluid like honey. When the particle reaches the
layer of fluid, it instantaneously sinks in it and the depth it reaches is
linked to the impact velocity. From then, the ball is glued to the
layer of fluid, the macroscopic contact begins, $q$ is set to
zero and $\gamma$ stores the impact velocity. As long as it is pushed, the
particle sinks deeper in the fluid and gets closer to the plane
($\gamma$ decreases). Then the particle is
pulled. From that moment on, it smoothly moves back from the fluid ($\gamma$ increases) and
comes unstuck from the layer of fluid when the
pulling forces have balanced the impact velocity and the pushing
forces ($\gamma$ reaches zero).
Note that, from~(\ref{PFD_contraint}), $\lambda$ can be interpreted as an
additional force, exerted by the plane on the particle, in order to
satisfy the constraint~(\ref{loi_choc}). It follows
from~(\ref{supp_lambda}) that the plane is allowed to act
on the particle through this force only if they are in macroscopic
contact.
\end{rmrk}

\begin{rmrk}[Radius]
\label{rmrk_radius}
This gluey particle model is built
in~\cite{Maury2007} as
the vanishing viscosity limit of the lubrication model~(\ref{DL_Flub_1}) where each
constant except the viscosity is taken equal to 1. Taking all constants into account leads to define $\gamma$ as the
limit of $\gamma_\mu=6\pi\mu ln(q)$ and 
the equation governing its evolution becomes
\begin{equation}
\label{evol_gamma_r}
\dot \gamma = -\frac{1}{r^2}\lambda.
\end{equation}
The larger $r$ is, the less the microscopic
distance $\gamma$ varies (the more it is difficult for the particle to
move). Note that, provided we are only interested in the macroscopic trajectory
$q$ of the
particle, the previous model (\ref{loi_choc})-(\ref{contrainte_q_gamma}) was valid for any radius: these
trajectories only depend on the sign of $\gamma$ (and not its value)
which is independent of $r$ from~(\ref{evol_gamma_r}).
\end{rmrk}

\begin{rmrk}[Viscous or not viscous ?]
Since this model is built by letting the viscosity go to zero, one may
wonder whether it models viscous fluids or not. To answer this
question, we consider the same experiment as in figure~\ref{schema_modele_limite} (pushing untill time 2 and
then pulling) for a particle falling on a plane coatted with different
viscous fluids. On figure~\ref{couche_visqueuse}, we
compare the trajectory given by the gluey particle model to the
trajectories computed for these systems where the viscous fluid layer
is modeled by~(\ref{PFD}).
\begin{figure}[hbtp]
\begin{center}
\resizebox{0.9\textwidth}{!}{\includegraphics{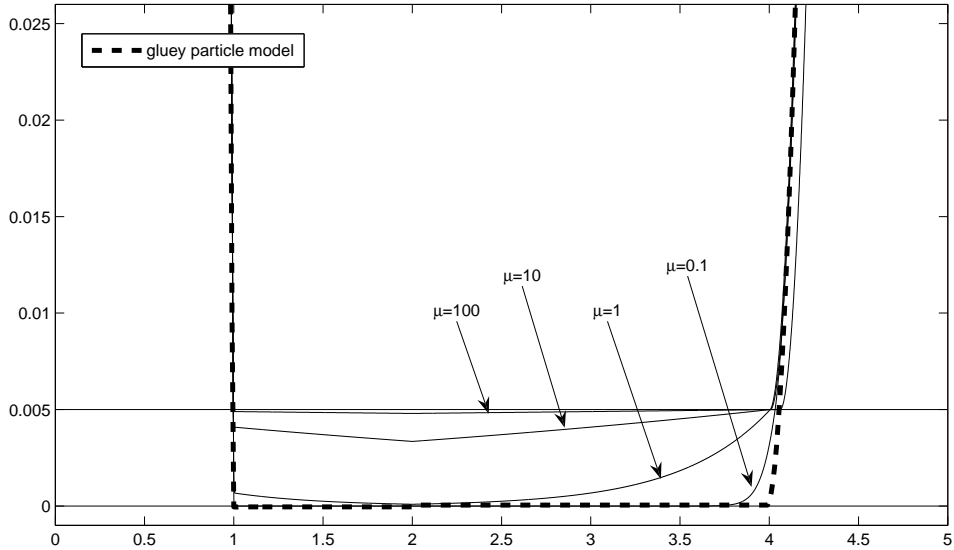}}
\end{center}
\caption{Comparison gluey particle model / layer of viscous fluid.}
\label{couche_visqueuse}
\end{figure}
Of course, trajectories converge to the limit model when the
viscosity goes to zero. We also observe that, from a macroscopic point of view, as long as we are interested in
hitting and unsticking times, the limit model seems to agree with all
trajectories. As a matter of fact, from this point of view, what is important is not the viscosity of the fluid but
whether it is viscous or not. However, the main difference between the
trajectories is the minimal distance reached by the
particle. Actually, small viscosities induce very small distances
and, consequently, the system can reach a domain wherein the initial
lubrication model~(\ref{PFD}) is no longer valid. To conclude, the gluey particle
model shall be employed to represent the
macroscopic behaviour of very viscous systems for which distances are
not too small.
\end{rmrk}

\subsection{Numerical algorithm}

We propose here an algorithm for
problem~(\ref{loi_choc})-(\ref{contrainte_q_gamma}). Let $h=T/N$ be the time step. The problem is initialized to $q^0>0$,
$u^0\in\IR$ and $\gamma^0=\lambda^0=0$. 
We denote by $q^n$,
$u^n$, $\gamma^n$ and $\lambda^n$
the computed values of $q$, $u$, $\gamma$ and $\lambda$ at time $t^n$.
We define $f^n$ by $f^n=\frac{1}{h}\int_{t^n}^{t^{n+1}}f(s)ds$. We have to
compute $q^{n+1}$, $u^{n+1}$, $\gamma^{n+1}$ and $\lambda^{n+1}$.

In order to compute $u^{n+1}$ and $\lambda^{n+1}$, we define the
discrete counterpart of $C_{q,\gamma}(t^n)$ the following way:
$$
\left|
\begin{array}{l}
\displaystyle K(q^n,\gamma^n)=\{v\virg q^n+hv\geq 0\} \text{ if }
\gamma^n=0,\vseq\\
\displaystyle K(q^n,\gamma^n)=\{v\virg q^n+hv= 0\}\text{ if }
\gamma^n<0.
\end{array}
\right.
$$
$K(q^{n},\gamma^n)$ is called the set of admissible velocities at
time $t^n$.
The collision law~(\ref{loi_choc}) and the Fundamental Principle of
Dynamics~(\ref{PFD_contraint}) then become,
$$
\left\{
\begin{array}{l}
\displaystyle u^{n+1/2}=u^n+hf^n, \vseq \\
\displaystyle u^{n+1}\in K(q^{n},\gamma^{n}) \virg \frac{1}{2}\left|u^{n+1}-u^{n+1/2}\right|^2_m=\min_{v\in
 K(q^{n},\gamma^{n})}\frac{1}{2}\left|v-u^{n+1/2}\right|^2_m,
\end{array}
\right.
$$
where $(v,w)_m=(mv,w)$. Note that
$u^{n+1/2}$ is the velocity the particle would have at time $t^{n+1}$
if there were no plane. $u^{n+1}$ is the projection of this a priori
velocity on the set of admissible velocities
$K(q^{n},\gamma^{n})$  for an adapted scalar product. From this
projection step, arises a Lagrange multiplier, denoted by $\lambda^{n+1}$ (positive
if $\gamma^{n}\geq 0$), and such that
$$m(u^{n+1}-u^{n+1/2})=h\lambda^{n+1}.$$
This can be rewritten as
\begin{equation}
\label{PFD_disc}
m\frac{u^{n+1}-u^n}{h}=mf^n+\lambda^{n+1},
\end{equation}
which is a discretization of~(\ref{PFD_contraint}).

Then, $\gamma^{n+1}$ is given by an explicit Euler discretization of~(\ref{evol_gamma}),
$$
\gamma^{n+1}=\gamma^n-h\lambda^{n+1}.
$$

This equation is valid while $\gamma^{n+1}$ is negative. If it becomes
strictly positive, it means that the particle has taken off at a time
$t^*\in]t_n,t_{n+1}[$. In that instance, $\gamma^{n+1}/m$ has integrated the
    force on  $]t^*,t^{n+1}[$ instead of $u^{n+1}$ which was fixed to
    zero. Therefore, in that case, we modify  $u^{n+1}$ and $\gamma^{n+1}$ the
    following way :
$$
\hbox{if }\gamma^{n+1}>0\virg u^{n+1}=\gamma^{n+1}/m \hbox{ and } \gamma^{n+1}=0.
$$
Finally the position $q^{n+1}$ is given by
$$
q^{n+1}=q^n+h u^{n+1}.
$$

\newpage
%
To sum up, the algorithm is the following :
\begin{lgrthm}[Particle/plane]
For all $n\geq 0$, let $q^n$, $u^n$, $\gamma^n$ and $\lambda^n$ be
given. We define $\displaystyle f^n=\frac{1}{h}\int_{t^n}^{t^{n+1}}f(s)ds$.
\begin{enumerate}
\item
\label{algo_u_apriori}
Computation of the a priori velocity, without taking the lubrication
force into account
$$\displaystyle u^{n+1/2}=u^n+hf^n.$$
\item
\label{algo_u}
Projection of the a priori velocity onto the set of {\it admissible velocities},
$$\displaystyle \bar u^{n+1}\in K(q^{n},\gamma^{n}) \virg
 \frac{1}{2}\left|\bar u^{n+1}-u^{n+1/2}\right|^2_m=\min_{v\in
 K(q^{n},\gamma^{n})}\frac{1}{2}\left|v-u^{n+1/2}\right|^2_m,
$$
$
\begin{array}{ll}
\hbox{where }& \displaystyle K(q,\gamma)=\{v\virg q+hv\geq 0\} \hbox{ if } \gamma=0, \vseq\\
 & \displaystyle K(q,\gamma)=\{v\virg q+hv= 0\} \hbox{ if } \gamma<0.
\end{array}
$

\vseq
From this projection step, we obtain $\lambda^{n+1}$.\vseq
\item 
\label{algo_gamma}
Updating of $\gamma$,
$$
\displaystyle \bar \gamma^{n+1}=\gamma^n-h\lambda^{n+1}.
$$
\item
\label{modif_decolle}
Modification if unsticking,
$$
\begin{array}{llll}
\text{if }\,\,\bar \gamma^{n+1}\leq0,
& u^{n+1}=\bar u^{n+1} & \text{and}
& \gamma^{n+1}=\bar \gamma^{n+1},\\
\text{if }\,\,\bar \gamma^{n+1}>0,
& u^{n+1}=\bar\gamma^{n+1}/m &\text{and}
& \gamma^{n+1}=0.
\end{array}
$$
\item
\label{algo_q}
Updating of $q$,
$$
q^{n+1}=q^n+hu^{n+1}.
$$
\end{enumerate}
\label{algo_cont_visc}
\end{lgrthm}


\begin{rmrk}[Coupling with fluid simulations]
\label{rmrk_coupling_algo}
This algorithm simulates collections of gluey
particles. Let us now suppose that the particles are embedded in a viscous
fluid. To make simulations taking the lubrication force into
account, a splitting method can be used to couple a fluid/particle
solver with the gluey particle algorithm. We denote by $\uu^n$ and
$\p^n$ the velocity and pressure fields into the fluid at time $t^n$.
Let $S$ be any fluid/particle
solver: from $\uu^n$, $q^{n}$ and $f^n$, $S$ computes the a priori
velocities of the particles, without taking the
lubrication force into account carefully. To couple the two algorithms we
propose to modify step~(\ref{algo_u_apriori}) of
algorithm~\ref{algo_cont_visc} writing:
$$
u^{n+1/2}=S(q^{n},\uu^n,f^n).
$$
\end{rmrk}

\subsection{Convergence result}

In this section, we establish a convergence result for the proposed
algorithm. To begin, we rewrite problem
(\ref{loi_choc})-(\ref{contrainte_q_gamma})
as
\begin{equation}
\left\{
\begin{array}{l}
\displaystyle m \dot q + \gamma = m\left(\dot q (0)+\int_0^t f(s)ds\right),\vseq\\
\displaystyle q\geq 0\virg \gamma\leq 0\virg q\gamma = 0, \vseq\\
\displaystyle q(0)=q^0>0\virg \dot q (0) =u^0 ,\vseq
\end{array}
\right.
\label{Pprime}
\end{equation}
which is formally equivalent
to the previous one (see~\cite{Maury2007}).

We recall that $h$ is the constant time step. We denote by
$q_h$ the piecewise affine function with $q_h(t^n)=q^n$. Similarly,
$\gamma_h$ is the piecewise affine function with $\gamma_h
(t^n)=\gamma^n$. We denote by $u_h$ the derivative of $q_h$, piecewise
constant equal to $u^{n+1}$ on $]t^n,t^{n+1}[$. Finally, we define
    $\lambda_h=-\dot \gamma_h$, piecewise constant. Note that, due to step~(\ref{modif_decolle}),
    $\lambda_h$ is generally not equal to $\lambda^{n+1}$ on
    $]t^n,t^{n+1}[$. We will denote by $\tilde\lambda^{n+1}=-(\gamma^{n+1}-\gamma^{n})/h$ its value on
    this interval. If the particle
    does not take off between times $t^n$ and $t^{n+1}$, no
    modification is made during step~(\ref{modif_decolle}) and we obtain
    $\tilde\lambda^{n+1}=\lambda^{n+1}$. The convergence theorem is the
    following:
\begin{thrm}
Let $f$ be integrable on $I=]0,T[$. When $h$ goes to zero, there
exists subsequences, still denoted by $(q_h)_h$, $(u_h)_h$, $(\lambda_h)_h$ and
$(\gamma_h)_h$, $q\in W^{1,1}(I)\cap {\mathcal C}(I)$ and $\gamma\in BV(I)$ such that
$$
\begin{array}{l}
u_h\longrightarrow u \hbox{ in } L^1(I), \vseq\\
q_h\longrightarrow q \hbox{ in } W^{1,1}(I) \hbox{ and } L^\infty(I) \hbox{
  with } \dot q =u, \vseq\\
\lambda_h\wsto\lambda \hbox{ in } {\mathcal M}(I), \vseq\\
\gamma_h\longrightarrow\gamma \hbox{ in } L^1(I) \hbox{ with } \dot\gamma=-\lambda,\vseq\\
\end{array}
$$
where $(q,\gamma)$ is solution to~(\ref{Pprime}).
\label{cv_schema}
\end{thrm}

\begin{rmrk}
Non-uniqueness for the limit problem (see~\cite{Maury2007} for
 counter-example) prevents from using the standard approach based on
 consistance and stability. Consequently, we use
 compactness methods and obtain convergence up to
 subsequences. However, in case $q$ has a finite number of zeros, the
 limit model admits a unique solution and therefore the convergence of
 the algorithm to problem~(\ref{Pprime}) is proved. Moreover, in that case, it can
be shown that (\ref{loi_choc})-(\ref{contrainte_q_gamma}) and (\ref{Pprime}) are
equivalent, in the sense that a solution to one of the problem is also
solution to the other (the demonstration of this result can be found
 in~\cite{these}). Consequently, under the a priori hypothesis that
 $q$ has a finite number of zeros, theorem~\ref{cv_schema} shows that
 algorithm~\ref{algo_cont_visc} converges to
 (\ref{loi_choc})-(\ref{contrainte_q_gamma}).
 For example, this hypothesis is verified if the external force $f$ changes
 of sign a finite number of times.
\end{rmrk}

\vseq
\noindent {\bf Proof of theorem~\ref{cv_schema}}
\vseq

To begin, note that a discrete form of the Fundamental Principle of
Dynamics~(\ref{PFD_contraint}) is verified:
\begin{equation}
\label{PFD_disc2}
\forall n\virg m\frac{u^{n+1}-u^n}{h}=mf^n+\tilde\lambda^{n+1}.
\end{equation}
Indeed, in case the particle does not take off between times $t^n$ and
$t^{n+1}$, the equality follows from~(\ref{PFD_disc}), together with
$\tilde\lambda^{n+1}=\lambda^{n+1}$. If the particle takes off, it
comes from~(\ref{PFD_disc}) and step~\ref{modif_decolle} of the
algorithm.

The proof will be devided into 4 steps.

\begin{enumerate}

\item {\bf Convergence of $q_h$ and $u_h$}

\begin{lmm}
\label{uh_bornee_Linf}
$(u_h)_h$ is bounded in  $L^\infty(I)$.
\end{lmm}

\noindent{\bf Proof :}
In case the particle does not take off, the projection
step~(\ref{algo_u}) gives
$$
|u^{n+1}|=|\bar u^{n+1}|\leq |u^{n+1/2}|\leq |u^n|+h|f^n|.
$$
In the other case,
it can be proved that
$
u^{n+1}\tilde\lambda^{n+1}\leq 0
$
and combining this with~(\ref{PFD_disc2}) gives the same
result. By summing up all these inequalities we obtain
$$
|u^{n+1}|\leq|u^0|+\int_0^T |f|,
$$
and the result follows from definition of $u_h$. \findemo

\begin{lmm}
\label{uh_bornee_BV}
$(u_h)_h$ is bounded in  $BV(I)$.
\end{lmm}

\noindent{\bf Proof :}
By lemma~\ref{uh_bornee_Linf}, the result will follow provided we
prove $\hbox{Var}(u_h)$
is bounded independently from $h$, where
$$\displaystyle \hbox{Var}(u_h)=\sum_{n=1}^{N-1}|u^{n+1}-u^{n}|.$$
To check this, we first split the sum and consider the sums between
indexes $p_1$ and $n_1$ where $t^{p_1}$ and $t^{n_1}$ are successive
unsticking times (See Fig.~\ref{lemme_uh_BV_notations}):
$$
\hbox{Var}_{[t^{p_1},t^{n_1}[}(u_h)=\sum_{n=p_1}^{n_1-1}|u^{n+1}-u^{n}|.
$$
The total variation of $u_h$ is made of a sum of such
terms.
\begin{figure}[hbtp]
\psfragscanon
\psfrag{p1}[l]{$t^{p_1}$}
\psfrag{n1}[l]{$t^{n_1}$}
\psfrag{p0}[l]{$t^{p_0}$}
\psfrag{n0}[l]{$t^{n_0}$}
\begin{center}
\resizebox{0.30\textwidth}{!}{\includegraphics{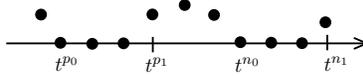}}
\end{center}
\caption{Proof of lemma~\ref{uh_bornee_BV} : notations.}
\label{lemme_uh_BV_notations}
\end{figure}

\noindent The idea behind the above decomposition is that, at each unsticking
time $t^{p_1}$, the velocity of the particle is small and that its
variations over $[t^{p_1},t^{n_1}[$ only depend on the integral of
    $f$ over the same interval. These terms will be summed up to obtain a bound
    on the total variation. 

More precisely, the bound for $\hbox{Var}_{[t^{p_1},t^{n_1}[}(u_h)$ can be found by analysing each jump
$|u^{n+1}-u^{n}|$, paying attention to what happens at time $t^n$
(hitting time, sticking time, unsticking time). For the sake of
    readibility, details of the
    computation are skipped here, they can be found in~\cite{these}. We find
$$
\hbox{Var}_{[t^{p_1},t^{n_1}[}(u_h)\leq
4\int_{t^{p_1-1}}^{t^{n_1}}|f(s)|ds.
$$
Summing up all these contributions and the bounding terms, we
obtain
$$
\hbox{Var}(u_h)\leq u^0+8\int_0^T |f(s)|ds,
$$
and  $\hbox{Var}(u_h)$
is bounded independently from $h$ as required.
\findemo

Lemma~\ref{uh_bornee_BV}, together with the compact embedding of
$BV(I)$ in $L^1(I)$ gives (up to a subsequence)
\begin{equation}
u_h\longrightarrow u \hbox{ in } L^1(I) \hbox{ with } u\in BV(I),
\label{u_cv_L1}
\end{equation}
$$
q_h\longrightarrow q \hbox{ in } W^{1,1}(I) \hbox{ with } \dot q=u.
$$
Uniform convergence of $q_h$ to $q$ then follows from the continuous embedding of $W^{1,1}(I)$ in
$L^\infty(I)$:
\begin{equation}
q_h\longrightarrow q \hbox{ in } L^\infty(I).
\label{q_cv_Linf}
\end{equation}
Finally, since $q_h$ is positive, we have $q\geq 0$ everywhere.\vseq

\item {\bf Convergence of $\gamma_h$}

\begin{lmm}
\label{lambdah_bornee_Lun}
$(\lambda_h)_h$ is bounded in $L^1(I)$.
\end{lmm}

\noindent{\bf Proof :}
By~(\ref{PFD_disc2}) and the fact that $\tilde\lambda^0=0$ we get
$$
\int_0^T|\lambda_h|\leq m \hbox{Var}(u_h)+\|f\|_{L^1(I)}.
$$
The result follows by combining this with lemma~\ref{uh_bornee_BV}.\findemo

By lemma~\ref{lambdah_bornee_Lun}, $(\lambda_h)_h$ is bounded in
${\mathcal M}(I)$, which implies that there exists a subsequence and
$\lambda\in {\mathcal M}(I)$ such that
$$
\lambda_h \wsto \lambda \hbox{ in } {\mathcal M}(I).
$$

Moreover, combining lemma~\ref{lambdah_bornee_Lun} with $\dot \gamma_h=-\lambda_h$ it comes
that $(\gamma_h)_h$ is bounded in $BV(I)$. This, together with compact
embedding of $BV(I)$ in $L^1(I)$, implies that there exists a subsequence
and $\gamma\in BV(I)$ such that
\begin{equation}
\gamma_h\longrightarrow\gamma \hbox{ in } L^1(I) \text{ and a.e.}
\label{gamma_cv_Lun}
\end{equation}
Since $\gamma_h$ is negative, it follows from this convergence result
that so is $\gamma$.
Finally, since $\dot \gamma_h=-\lambda_h$, we can check that $\dot \gamma=-\lambda$ in ${\mathcal M}(I)$.\vseq

\item {\bf Continuous FPD}

We are now going to prove that $\displaystyle m \dot q + \gamma = m\left(\dot q
(0)+\int_0^t f(s)ds\right)$ almost everywhere on $I$. 

In order to do so, the first step is to
prove that~(\ref{PFD_contraint}) is verified in the sense of
distributions. From~(\ref{PFD_disc2}) it follows that
\begin{equation}
\forall \varphi\in {\mathcal D}(I)\virg \left<m\dot u_h,\varphi \right>=\sum_{n=1}^{N-1}mhf^{n}\varphi(t^n)+\sum_{n=1}^{N-1}h\tilde\lambda^{n+1}\varphi(t^n).
\label{PFD_disc3}
\end{equation}
We are going to pass to the limit in this equation. By~(\ref{u_cv_L1}), $\left<m\dot u_h,\varphi \right>$
converges to $\left<m\dot u,\varphi \right>$. To study the first term
of the right-hand side, we write
$$
\displaystyle h\sum_{n=1}^{N-1}f^{n}\varphi(t^n)
=
\int_0^T f(s)\varphi(s)ds +
\sum_{n=1}^{N-1}\int_{t^{n}}^{t^{n+1}}f(s)\left[\varphi(t^n)-\varphi(s)\right]ds
-\int_{t^{0}}^{t^1}f(s)\varphi(s)ds.
$$
The convergence to zero of the sum over $n$ comes from uniform
continuity of $\varphi$. Combining this with $|t^1-t^{0}|=h$ gives
$$
\sum_{n=1}^{N-1}mhf^{n}\varphi(t^n) \longrightarrow m\int_0^T f(s)\varphi(s)ds\hbox{ when } h\rightarrow 0.
$$
The argument for the last term is similar. We write
$$
\displaystyle \sum_{n=1}^{N-1}h\tilde\lambda^{n+1}\varphi(t^n)
 =\int_0^T \lambda_h(s)\varphi(s)ds +
\sum_{n=1}^{N-1}\int_{t^n}^{t^{n+1}}\lambda_h(s)\left[\varphi(t^n)-\varphi(s)\right]ds
-\int_{t^0}^{t^1}\lambda_h(s)\varphi(s)ds.
$$
The convergence to zero of the sum over $n$ comes from uniform
continuity of $\varphi$ and lemma~\ref{lambdah_bornee_Lun}, and the last
term is equal to zero for all $h$. This, together with
lemma~\ref{lambdah_bornee_Lun} gives
$$
\sum_{n=1}^{N-1}h\tilde\lambda^{n+1}\varphi(t^n)\longrightarrow
\langle\lambda,\varphi\rangle=-\langle\dot\gamma,\varphi\rangle \hbox{ when } h\rightarrow 0.
$$
Finally, passing to the limit in~(\ref{PFD_disc3}) we obtain
$$
\langle m\ddot q-\dot\gamma,\varphi\rangle=\langle mf,\varphi\rangle\virg
\forall \varphi\in {\mathcal D}(I),
$$
as required.

Then, by density of ${\mathcal D}(I)$ in ${\mathcal C}^0_0(I)$ and the fact
that $m\dot q - \gamma$ is in $BV(I)$, we get
$$
m\ddot q-\dot\gamma=mf \hbox{ in } {\mathcal M}(I).
$$

Integrating this equality over $[0,t[$ (Stieltjes integral of $BV$
    functions) we obtain
$$
(m\dot q-\gamma)(t^+)-(m\dot q-\gamma)(0^-)=\int_0^t mf,
$$
and the result follows from this, by using $\gamma(0^-)=0$ and
a.e. continuity of $m\dot q-\gamma$.\vseq

\item {\bf Proof of $q\gamma=0$}

To prove that $(q,\gamma)$ is solution to (\ref{Pprime}),
it remains to show that $q\gamma=0$ almost everywhere. For all $n$ we
have $q^n\gamma^n=0$. However, $q_h\gamma_h$ is not identically equal
to zero. We build new
functions $\tilde q_h$ and $\tilde\gamma_h$, piecewise constant, with
respective values $q^n$ and $\gamma^n$ on $]t^n,t^{n+1}[$. We now have
    $\tilde q_h\tilde\gamma_h=0$ and simple
    computations give
$$
\|\tilde q_h-q\|_{L^\infty(I)}\leq\|\tilde
q_h-q_h\|_{L^\infty(I)}+\|q_h-q\|_{L^\infty(I)}
\leq h\|u_h\|_{L^\infty(I)} + \|q_h-q\|_{L^\infty(I)}
$$
and
$$
\|\tilde \gamma_h-\gamma\|_{L^1(I)}\leq\|\tilde
\gamma_h-\gamma_h\|_{L^1(I)}+\|\gamma_h-\gamma\|_{L^1(I)}
\leq \frac{h}{2}\|\lambda_h\|_{L^1(I)} + \|\gamma_h-\gamma\|_{L^1(I)}.
$$
Combining the first inequality with lemma~\ref{uh_bornee_Linf}
and~(\ref{q_cv_Linf}) gives uniform convergence of $\tilde q_h$ to
$q$. Putting together the second inequality,
lemma~\ref{lambdah_bornee_Lun} and~(\ref{gamma_cv_Lun}), we see that
$\tilde\gamma_h$ converges to $\gamma$ in $L^1(I)$ which implies that
the sequence converges up to a subsequence almost everywhere on $I$. Finally, letting $h$
go to zero in $\tilde q_h\tilde\gamma_h=0$ gives $q\gamma=0$ almost
everywhere as required.
\end{enumerate}

This completes the proof of theorem~\ref{cv_schema}. \findemo

\subsection{Validation: coupling with a
  fluid/particle solver}
\label{sec:couplage}

We consider the same
experiment that the one considered in section~\ref{sec:plane_part},
  without inertia. The radius of the particle is taken equal to $1$ and the
  viscosity of the fluid is $\mu=3$. The balance of forces reads
\begin{equation}
\label{eq_forces}
\forall t\virg F_{lub}(q(t))+f(t)=0,
\end{equation}
where $f(t)=-2$ until time 2 and $f(t)=2$ if $t>2$. 

\begin{figure}[hbtp]
\centering
\includegraphics[width=0.49\textwidth]{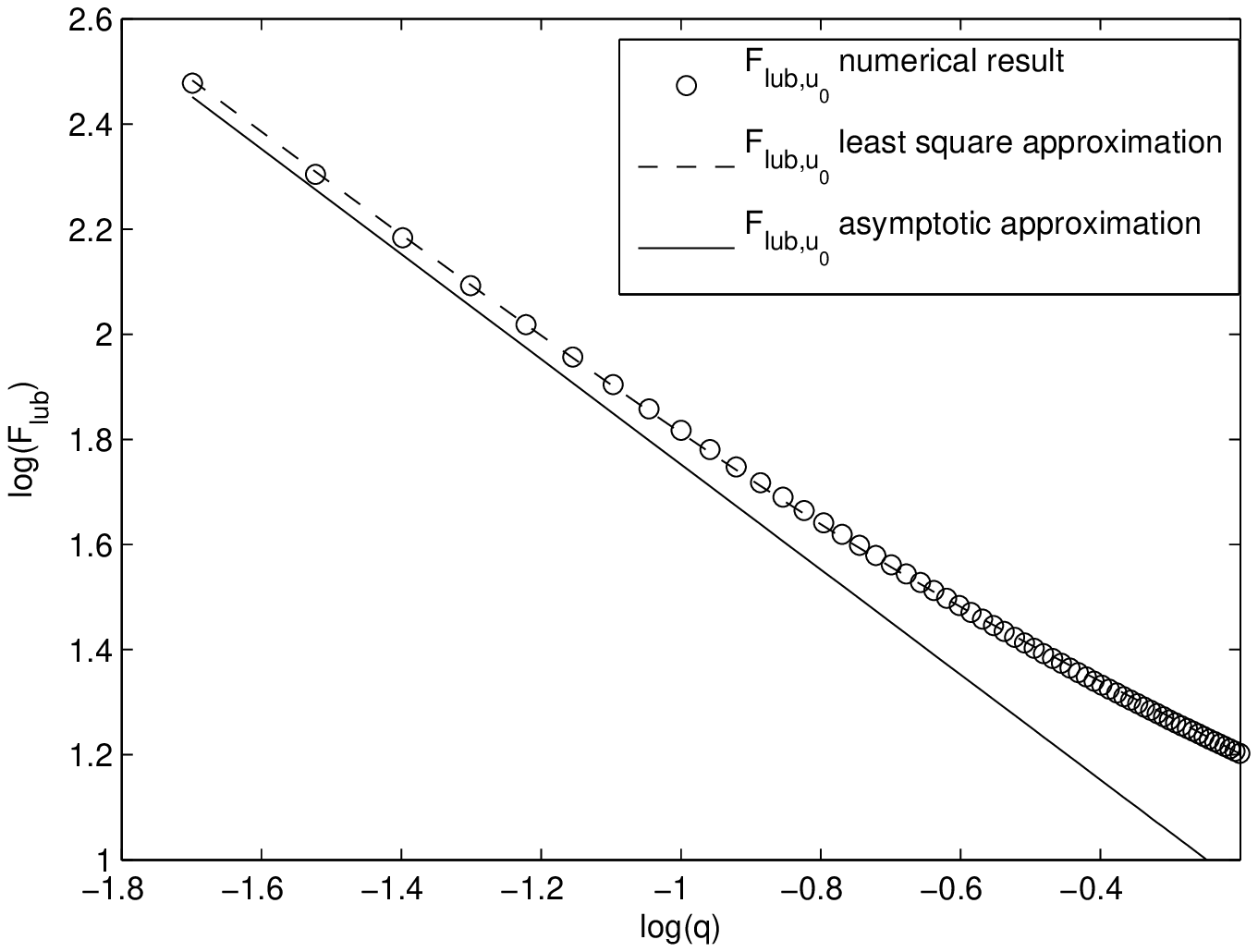}
\includegraphics[width=0.49\textwidth]{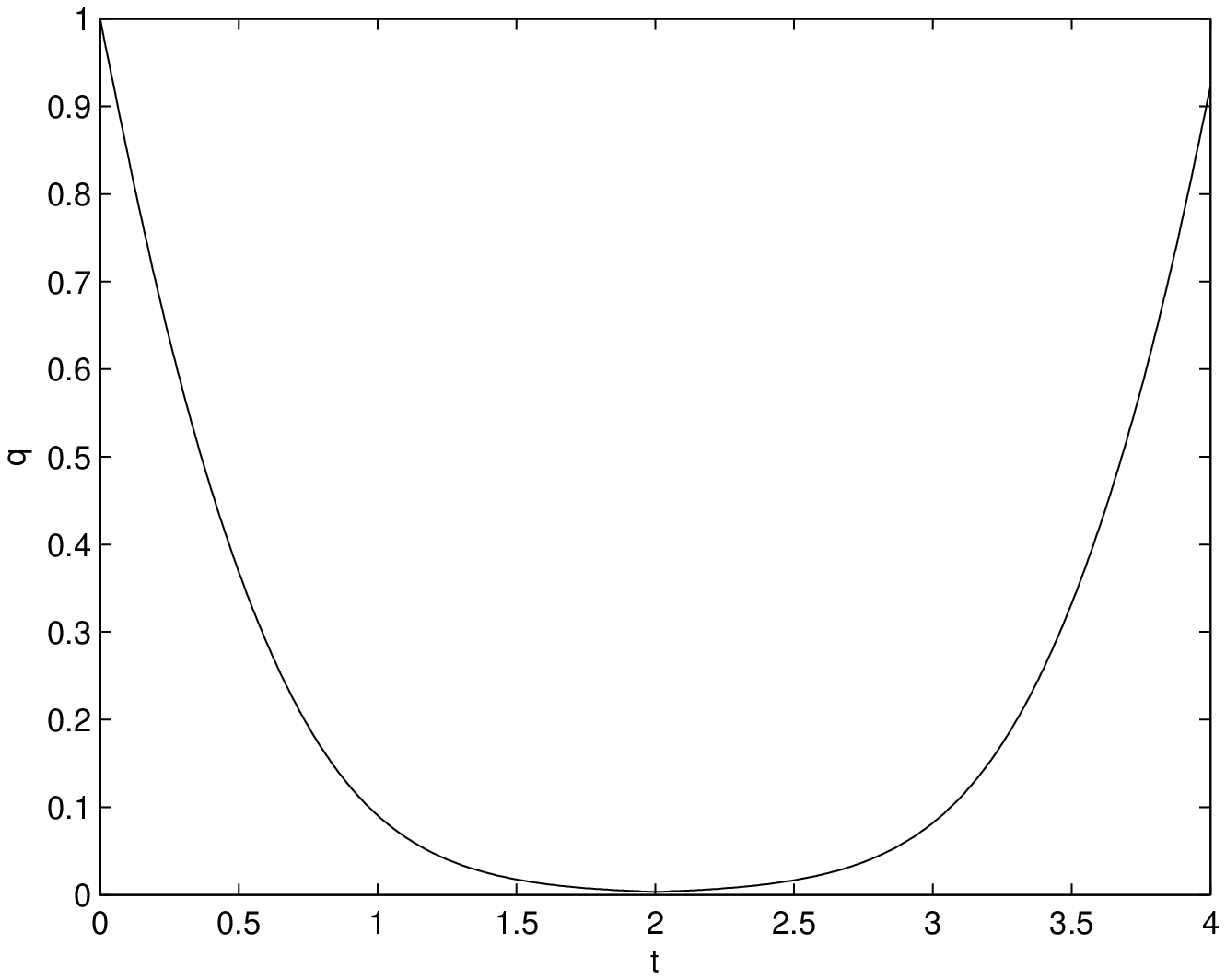}
\caption{Approximation of the lubrication force (left) and reference
  solution (right).}
\label{reference}
\end{figure}

To obtain a reference solution we first compute, as accurately as
possible, the map $q\rightarrow F_{lub,u_0}(q)$ for a given velocity
$u_0=-1$ and $q\in [0,1]$. To do so, we begin with computing $F_{lub,u_0}(q_k)$ where
$(q_k)_{k=1..M}$ is a regular subdivision of interval $[0,1]$. This is
done, for each $q_k$, solving the Stokes problem in the fluid with Dirichlet
boundary conditions and computing the force $F_{lub,u_0}(q_k)$ exerted by the fluid on
the particle. The computations are carried out in tree-dimensions using an
axisymmetric formulation and the Finite-Element solver {\tt
  FreeFem++}. On the left side of figure~\ref{reference}, we plot the numerical
results obtained (circles). They agree with the asymptotic
expansion~(\ref{DL_Flub_1}) for small distances (solid line). Finally, the map
$q\rightarrow F_{lub}(q)$ is approximated using a least square approximation
of the numerical results by a polynomial of degree $3$ (dashed line).

The reference solution is obtained discretizing the time interval and
computing the velocity of the particle $u^n$ at each time-step. We write that $u^n=\alpha^n
u_0$  and, using the linearity of the
lubrication force with respect to the velocity, we compute $\alpha^n$ as
the solution to 
$$
\alpha^n F_{lub,u_0}(q^n)+f(t^n)=0.
$$
The trajectory obtained is plotted against time on the right side of figure~\ref{reference}.

\begin{figure}[hbtp]
\centering
\includegraphics[width=0.5\textwidth]{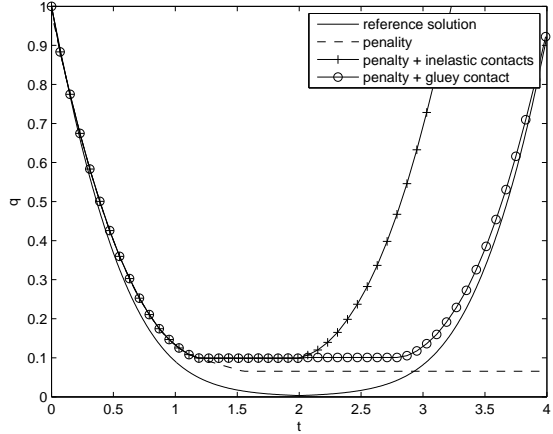}
\caption{Comparison of the numerical solutions for different contact models.}
\label{compare_meth}
\end{figure}

We now want to observe the influence of the method employed to deal
with contacts in fluid/particle simulations. To do so, we use an
axisymmetrical version of the fluid/particle solver implemented with {\tt
  FreeFem++} and described
in~\cite{GTN}. On figure~\ref{compare_meth}, we plot the solution
given by this solver for a mesh size $\delta x = r/10$ (dashed line). We can observe
that the particle remains glued. Indeed, due to the space
discretization, the characteristic function representing the rigid
particle ends up with touching the boundary of the domain and the
Dirichlet boundary condition prevents it from
taking off. Consequently, it is necessary to deal with the problem of
contact and to prevent the characterictic function from intersecting the
boundary of the domain. Two methods are tested: the fluid/particle
solver is coupled with an inelastic contact algorithm and with the
gluey contact model. The coupling is performed using the splitting
strategy described in remark~\ref{rmrk_coupling_algo}. In each case, the constraint for the distance is set to
$q\geq \eta$ with $\eta=\delta x$. The numerical results are compared
on figure~\ref{compare_meth}. We observe that, for the inelastic model
(solid line with crosses), the
particle takes off as soon as it is pulled. To the contrary, using the
gluey contact model (solid line with circles), the particle remains glued and the trajectory
finally joins up with the reference one. This is a validation of the gluey
particle model and it emphasizes the necessity to take the lubrication force
into account when dealing with contacts.

\begin{figure}[hbtp]
\centering
\includegraphics[width=0.49\textwidth]{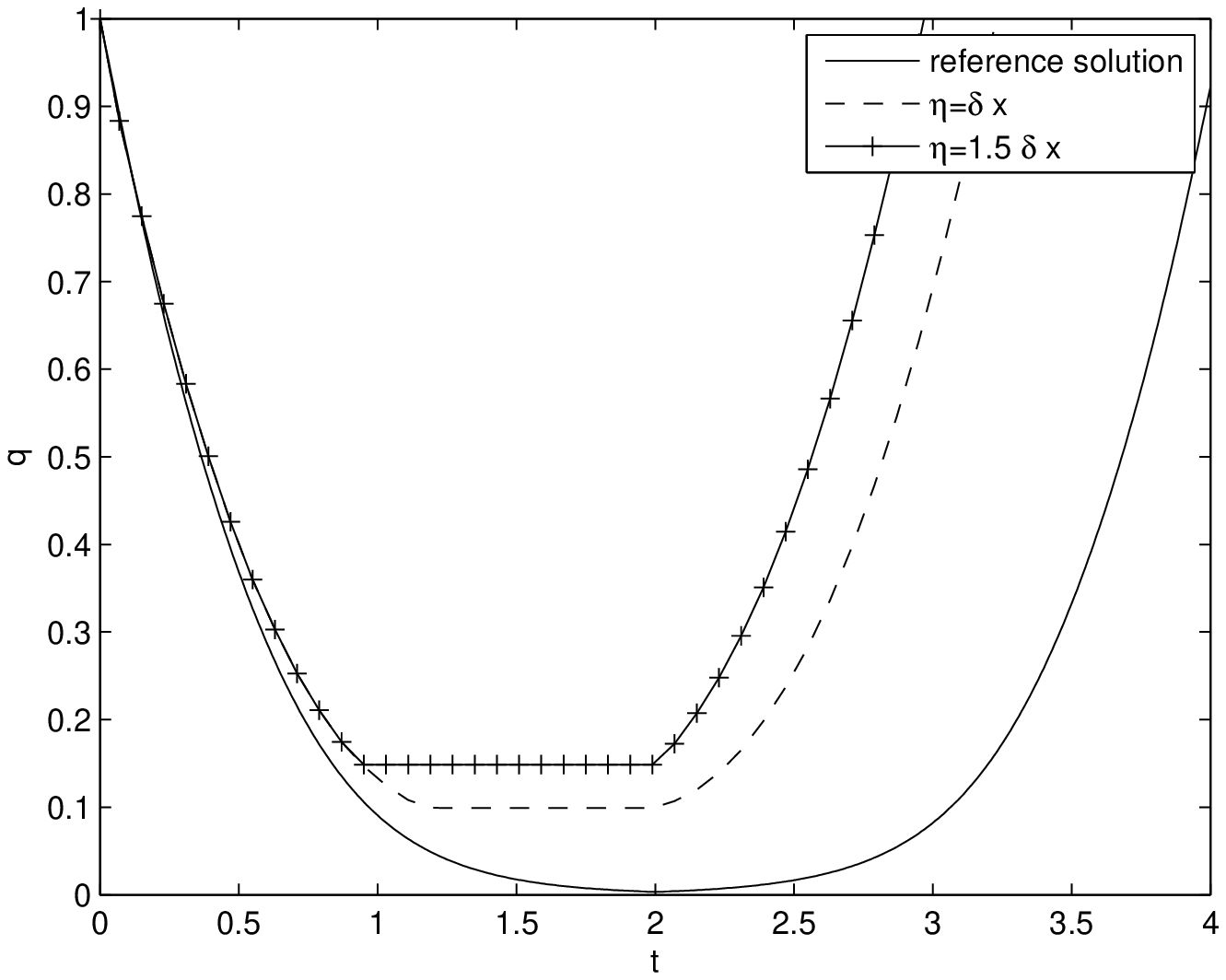}
\includegraphics[width=0.49\textwidth]{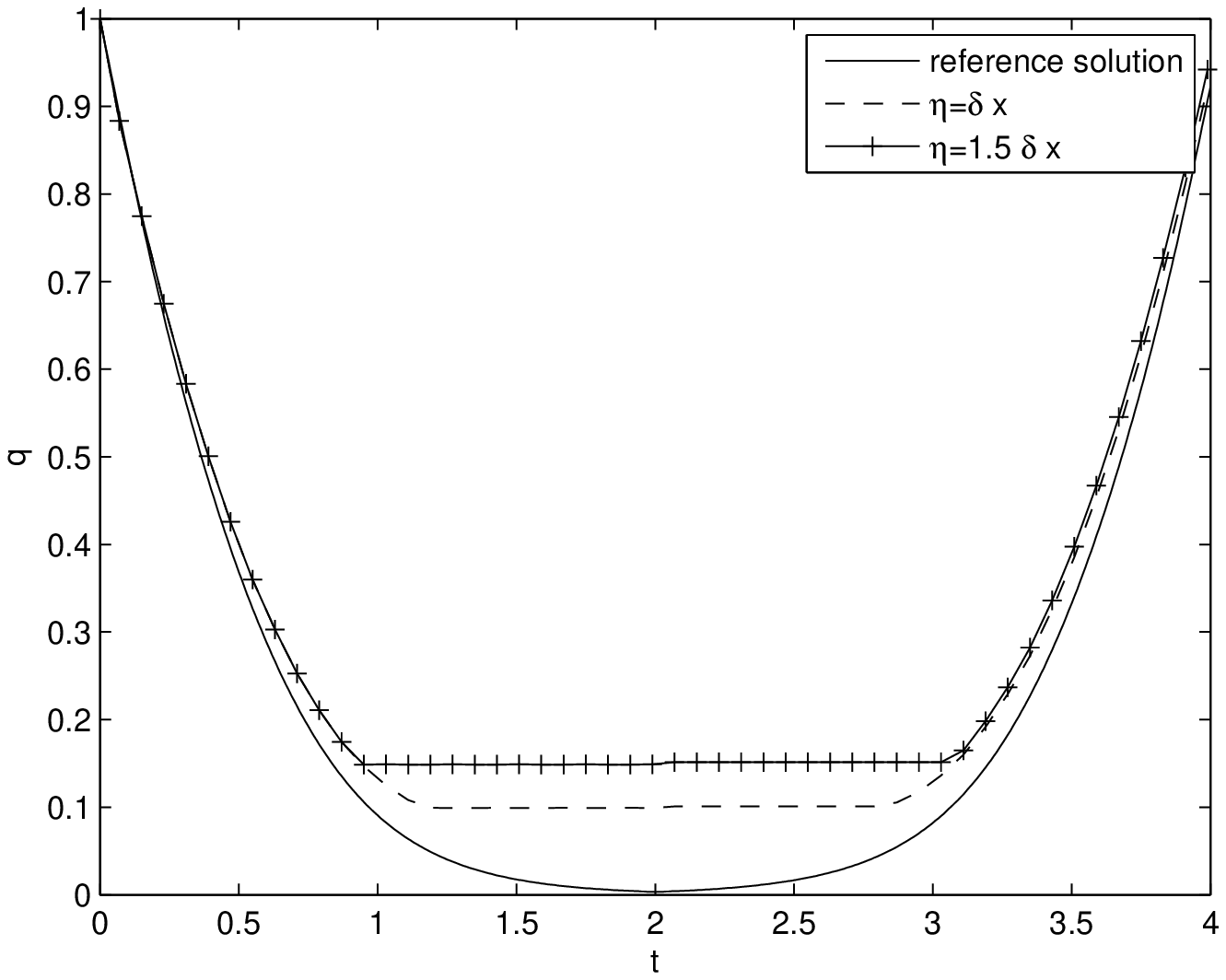}
\caption{Impact of $\eta$ on the numerical solution for inelastic
  (left) and gluey (right) contact.}
\label{influence_seuil}
\end{figure}

Finally, we observe on figure~\ref{influence_seuil} the behaviour of
the two contact models with respect to the parameter $\eta$ which is the
minimal distance allowed between the particle and the plane. We can see that the trajectories obtained
for different $\eta$ separates after unsticking time when using the inelastic contact
model (left side of the figure). This is due to the fact that, for this model, the particle unsticks as soon as
it is pulled. To the contrary, the gluey particle model
is not so sensible to parameter $\eta$ (right side of the figure).

\subsection{Extension to rough solid surfaces}
\label{sec:rough_solids}

\begin{figure}[hbtp]
\psfragscanon
\psfrag{r}[l]{$r$}
\psfrag{d}[l]{$q$}
\psfrag{b}[l]{$q_s<r_s$}
\psfrag{bard}[l]{$q_s<r_s$}
\psfrag{rs}[l]{$r_s$}
\psfrag{rsmbard}[l]{$r_s-q_s$}
\begin{center}
\resizebox{0.40\textwidth}{!}{\includegraphics{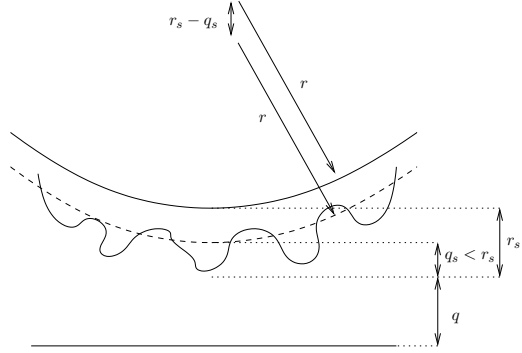}}
\end{center}
\caption{Equivalent smooth sphere.}
\label{sphere_equiv}
\end{figure}

As already said, it has been proved that smooth solids can not undergo contact. However, from our
experience, we know that the particle should touch the
plane in finite time. One of the reasons explaining this behaviour is that physical particles are not
smooth. Recent
experiments described in~\cite{Vinogradova2006} show
that the lubrication force exerted on a rough particle is the one
would be exerted on a shifted smooth particle:
\begin{equation}
\label{F_lub_rug}
F_{lub,rough}\sim -6\pi\mu r^2\frac{\VV}{q+q_s},
\end{equation}
where $q_s<r_s$ (see Fig.~\ref{sphere_equiv}).
Due to the lubrication force, this equivalent smooth sphere can not undergo
contact with the plane ($q+q_s$ doesn't go to zero in finite
time) but the real surfaces can collide ($q$ can go to zero). 

\begin{figure}[hbtp]
\psfragscanon
\psfrag{q}[l]{ $q$}
\psfrag{ra}[l]{ $r_{s,1}$}
\psfrag{rb}[l]{ $r_{s,2}$}
\centering
\includegraphics[width=0.4\textwidth]{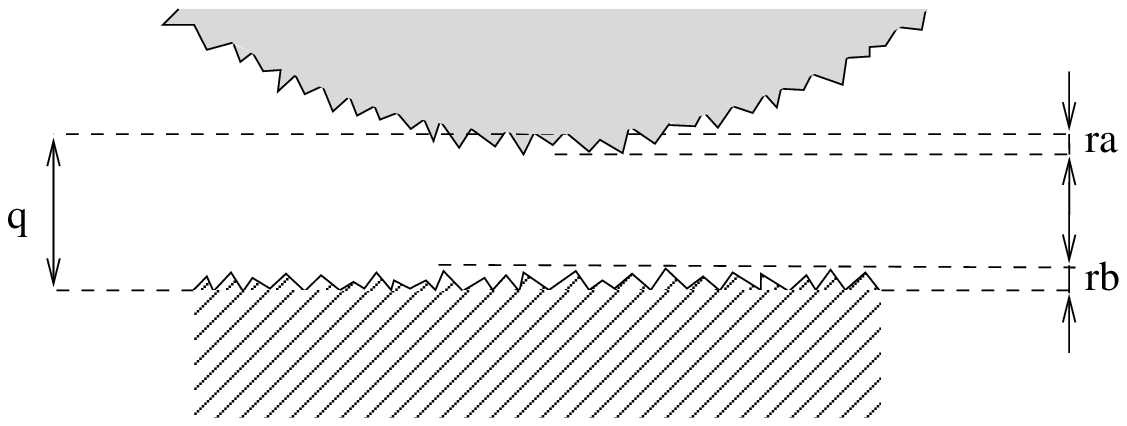}
\qquad
\includegraphics[width=0.5\textwidth]{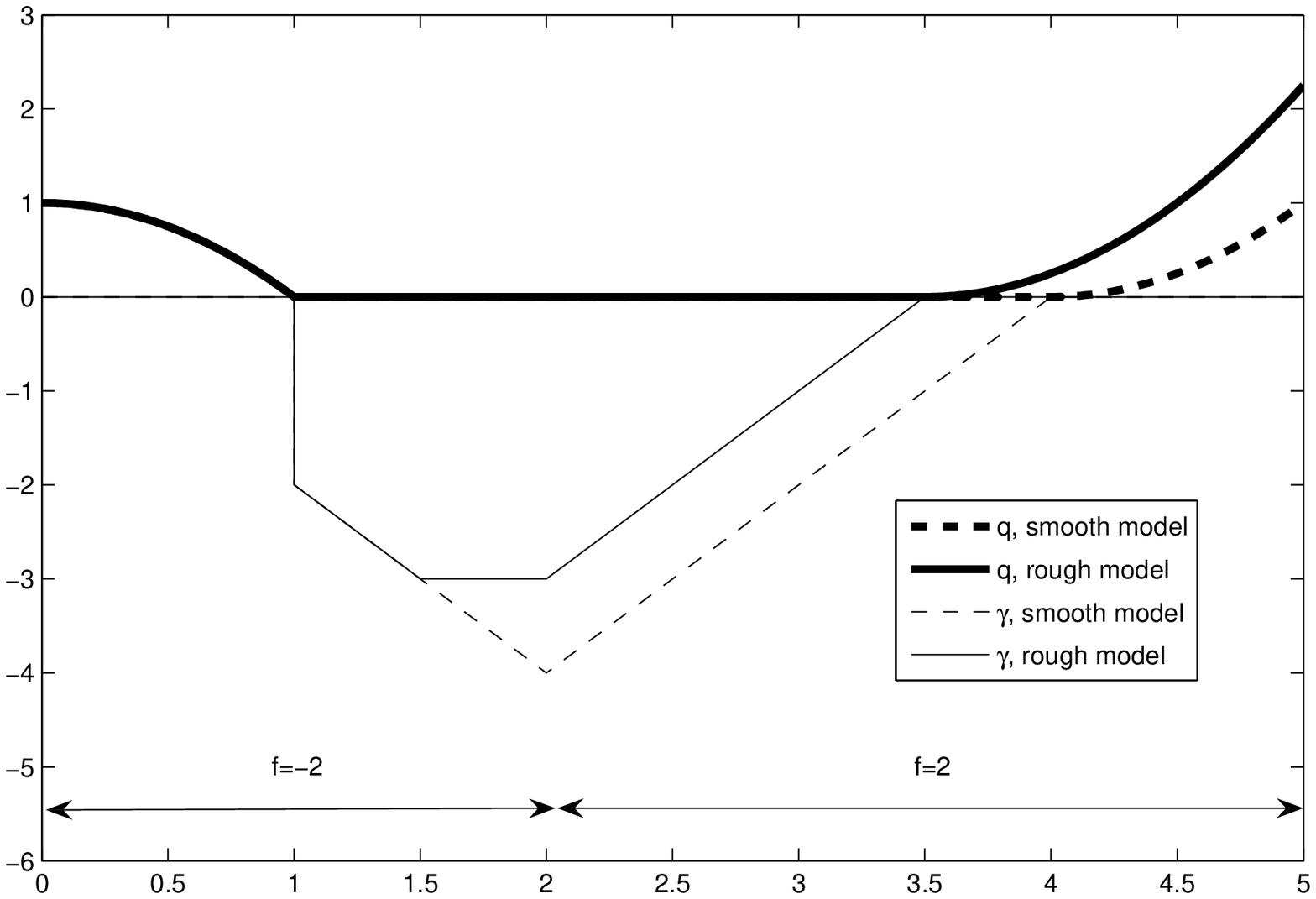}
\caption{Rough solids: notations (left) and gluey particle model (right).}
\label{rugueux}
\end{figure}

Taking these results into account in the gluey particle model, we consider
that, as soon as $q=r_{1,s}+r_{2,s}$ (see
notations on figure~\ref{rugueux}), there exists a real solid/solid
contact. During this contact, the forces
acting on the particle are not registered
anymore. To model such a behaviour, it suffices
to recall that $\gamma$ is the limit of $\gamma_\mu=6\pi\mu ln(q_\mu)$
(see remark~\ref{rmrk_radius}) and to impose 
$$\gamma\geq
6\pi\mu ln(r_{s,1}+r_{s,2}).
$$
The trajectory computed for this model is
plot on figure~\ref{rugueux}. We can observe that the rough particle
takes off before the smooth one.
Contrary to what has been said in remark~\ref{rmrk_radius} for the
smooth case, it is now important to know the value of $\gamma$ in order to
truncate it. Therefore, it is essential to take the radius into
account in its evolution and to use equation~(\ref{evol_gamma_r}) in
the gluey particle model:
$$
\dot \gamma = -\frac{1}{r^2}\lambda.
$$
In that case, the trajectory of the particle depends on $r$.

From an algorithmic point of view, this rough model can easily
be taken into account in algorithm~\ref{algo_cont_visc} by changing step~(\ref{algo_gamma}) in
$$
\displaystyle \gamma^{n+1}=\gamma^n-\frac{h}{r^2}\lambda^{n+1},
$$
where $r$ is the radius of the sphere and by adding the
following~(\ref{modif_decolle}b) step:
$$
\text{if}\,\,\, \gamma^{n+1}< \gamma_{min}\virg\gamma^{n+1} = \gamma_{min}.
$$

\section{Multi-particle case}

\subsection{Modelling}

\begin{figure}[hbtp]
\psfragscanon
\psfrag{qi}[l]{$\xx_i$}
\psfrag{qj}[l]{ $\xx_j$}
\psfrag{eij}[l]{$\ee_{ij}(\xx)$}
\psfrag{Dij}[l]{$D_{ij}(\xx)$}
\begin{center}
\resizebox{0.25\textwidth}{!}{\includegraphics{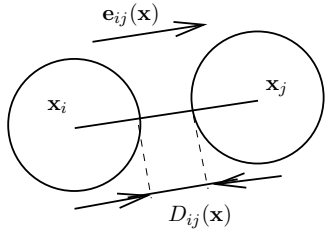}}
\end{center}
\caption{Particles i et j : notations.}
\label{schema_dim2_ij}
\end{figure}

We generalize the gluey particle model
(\ref{loi_choc})-(\ref{contrainte_q_gamma}) to the multi-particle
case. We consider a system of $N$ spherical particles in three-dimensions. $\xx_i$ stands
for the position of the center of particle $i$ in $\IR^3$ and  $\ff_i\in\IR^3$
for the external force exerted on it. 
Let $\xx\in\IR^{3N}$ be defined by $\xx=(\ldots,\xx_i,\ldots)$ and
$\ff\in\IR^{3N}$ by $\ff=(\ldots,\ff_i,\ldots)$. We denote by $D_{ij}$ the signed distance between
particles $i$ and $j$, and $\ee_{ij}$ by $\ee_{ij}(\xx)=(\xx_j-\xx_i)/\|\xx_j-\xx_i\|$
(see Fig.~\ref{schema_dim2_ij}). We define $M$ as the
mass matrix of dimension $3N\times 3N$,
$M=diag(\ldots,m_i,m_i,m_i,\ldots)$. Vector $\GGij\in\IR^{3N}$ is the gradient of
distance $D_{ij}$ with respect to the positions of the particles:
$$
\begin{array}{cccccc}
\GGij(\xx)=\nabla_\xx D_{ij}(\xx)=&(\ldots ,0,&-\ee_{ij}(\xx)&, 0,\ldots,0,&\ee_{ij}(\xx)&, 0,\ldots,0)^t.\\
 & &i& &j&\\
\end{array}
$$
\vspace{-1cm}

In that context, there are $N(N-1)/2$ pair of particles and we denote
by $\ggamma=(\ldots,\gammaij,\ldots)\in\IR^{N(N-1)/2}$ the associated
sticking variables: $\gammaij$ is stricly negative if particles $i$ and $j$ are
glued. Then, using the fact that $\displaystyle
\frac{dD_{ij}(\xx)}{dt}=\GGij(\xx)\cdot\dot\xx$, we define the
following space of admissible velocities: 
$$
C_{\xx,\ggamma}(t)=\left\{
\VV \in\IR^{3N}\hbox{ s.t.}
\left|
\begin{array}{l}
\displaystyle \GGij(\xx)\cdot\VV=0 \hbox{ if }  \, \gammaij(t^-)<0 \vseq\\
\displaystyle \GGij(\xx)\cdot\VV\geq 0 \hbox{ if } \, \gammaij(t^-)=0\hbox{, } D_{ij}(t)=0
\end{array}
\right.
\right\}.
$$
To finish with notations, we denote by $\llambda=(\ldots,\lambdaij,\ldots)\in\IR^{N(N-1)/2}$ the
vector made of the Lagrange multipliers associated to these $N(N-1)/2$
constraints. 

The multi-particle model is the natural counterpart of the
particle/plane one:
\begin{equation}
\label{model_N}
\left\{
\begin{array}{l}
\displaystyle \xx\in (W^{1,\infty}(I))^{3N}\virg \dot \xx\in (BV(I))^{3N}\virg\ggamma\in
(BV(I))^{N(N-1)/2}\virg\llambda\in ({\mathcal
  M}(I))^{N(N-1)/2}\hbox{,}\vseq\\
\displaystyle \dot \xx(t^+) = P_{C_{\xx,\ggamma}(t)}\dot \xx(t^-),\vseq \\
\displaystyle M \ddot \xx =  M \ff + \sum_{i<j}\lambdaij \GGij(\xx),\\
\displaystyle \hbox{supp}(\lambdaij)\subset \{t\virg D_{ij}(t)=0\}\hbox{ for all } i,j,\vseq\\
\displaystyle \dot \ggamma = -\llambda, \vseq\\
\displaystyle D_{ij}\geq 0\virg\gammaij\leq 0\hbox{ for all } i,j,\vseq\\
\displaystyle \xx(0)=\xx^0 \hbox{ st. } D_{ij}(0)>0\hbox{ for all } i,j\virg
\dot \xx (0) =\uu^0 \virg \ggamma(0)=0_{\IR^{N(N-1)/2}}.
\end{array}
\right.
\end{equation}

\begin{rmrk}
The Lagrange multiplier $\lambdaij$, associated to the constraint
between particles $i$ and $j$, is activated (non zero) only if these particles
are in contact. The additional force due to this contact is
equal to $\lambdaij \GGij(\xx)$. From the expression of $\GGij(\xx)$, we get
that this force only concerns the particles involved in the contact: it is equal to $-\lambdaij\ee_{ij}(\xx)$ on particle $i$ and $\lambdaij\ee_{ij}(\xx)$
on particle $j$.
\end{rmrk}
\begin{rmrk}[Roughness and radius]
As for
the particle/plane case (see section~\ref{sec:rough_solids}), roughness can be taken
into account by imposing a threshold on $\ggamma$:
$$
6\pi\mu ln(r_{i,s}+r_{j,s})\leq\gammaij \text{ for all } i,j,
$$
where $r_{l,s}$ is the size of roughness of particle $l$.
As noticed in the particle/plane case, it is now important to
take the radius of the particles into account in the evolution of
$\ggamma$. To do so, in the same way as in the particle/plane case, we come back to the way the gluey
particle model has been built and take into account all the
constants involved in the first order assymptotic
developpement of the lubrication force exerted between two particles~(\ref{DL_Flub_2}).
We obtain the following evolution equation for $\ggamma$:
$$
\dot \ggamma = -R\llambda,
$$
where $R$ is
the diagonal matrix of dimension $N(N-1)/2$ with coefficients $R_{ij,ij}=(r_i+r_j)^2/(r_i^2r_j^2)$.
\end{rmrk}

\subsection{Algorithm}

Let $h$ be the time step. We denote by $\VV^n=(\ldots,
\VV_i^n,\ldots)\in\IR^{3N}$ the approximated velocities of the
particles at time $t^n=nh$. Let $\xx^n$, $\ggamma^n$ and $\llambda^n$
be the respective approximations of $\xx$, $\ggamma$ and
$\llambda$ at time $t^n$.

The
discretization of the continuous constraints $C_{\xx,\ggamma}(t^n)$ is
inspired by~\cite{Maury2006} and corresponds to a first order
approximation of the constraints:
$$
K(\xx^{n},\ggamma^{n})=\left\{
\VV\in\IR^{3N} \text{ s.t.}
\left|
\begin{array}{l}
D_{ij}(\xx^n)+h\GGij(\xx^n)\cdot\VV\geq 0 \hbox{ if } \gammaij^n = 0\vseq\\
D_{ij}(\xx^n)+h\GGij(\xx^n)\cdot\VV = 0 \hbox{ if } \gammaij^n < 0\\
\end{array}
\right.
\right\}.
$$

Using this discrete space of admissible velocities, the time
discretization of~(\ref{model_N}) is now a direct adaptation of
algorithm~\ref{algo_cont_visc} to the multi-particle case. 
\newpage
\begin{lgrthm}[Multi-particle]
For all $n\geq 0$, let $\xx^n$, $\VV^n$, $\ggamma^n$ and $\llambda^n$ be
given. We define $\displaystyle \ff^n=\frac{1}{h}\int_{t^n}^{t^{n+1}}\ff(s)ds$.
\begin{enumerate}
\item
\label{algo_N_u_apriori}
Computation of the a priori velocity, without taking the lubrication
force into account
$$\displaystyle \VV^{n+1/2}=\VV^n+h\ff^n.$$
\item
\label{algo_N_u}
Projection of the a priori velocity on the set of {\it admissible
  velocities},
$$\displaystyle \VV^{n+1}\in K(\xx^{n},\ggamma^{n}) \virg \frac{1}{2}\left|\VV^{n+1}-\VV^{n+1/2}\right|^2_M=\min_{\VV\in
 K(\xx^{n},\ggamma^{n})}\frac{1}{2}\left|\VV-\VV^{n+1/2}\right|^2_M.
$$
From this projection step, we obtain $\llambda^{n+1}$.\vseq
\item 
\label{algo_N_gamma}
Updating of $\ggamma$,
$$\begin{array}{l}
\displaystyle \ggamma^{n+1}=\ggamma^n-h\llambda^{n+1},\vseq\\
\displaystyle \hbox{if } \,\,\gammaij^{n+1}>0\virg \gammaij^{n+1}=0.\\
\end{array}
$$
\item
\label{algo_N_q}
Updating of $\xx$,
$$
\xx^{n+1}=\xx^n+h\VV^{n+1}.
$$
\end{enumerate}
\label{algo_cont_visc_N}
\end{lgrthm}

\begin{rmrk}
In the same way as in section~\ref{sec:rough_solids}
and remark~\ref{rmrk_coupling_algo} for the particle/plane case, this algorithm can be
extended to rough solids and coupled with fluid/particle solvers using
a splitting strategy.
\end{rmrk}

\begin{rmrk}[Obstacles]
Suppose there exists $N_0$ fixed obstacles (walls of a box containing
the particles for example). It is straightforward to add the $NN_0$
new constraints in $K(\xx^{n},\ggamma^{n})$. Now, suppose these
obstacles are moving with a prescribed velocity. We denote by $\yy^{n+1}$
the (known) vector giving their position at time $t^{n+1}$. The space of
admisible velocities becomes:
$$
K(\xx^{n},\yy^{n+1},\ggamma^{n})=\left\{
\VV\in\IR^{3N} \text{ s.t.}
\left|
\begin{array}{l}
\hbox{Pairs }(i,j)\hbox{ particle/particle :}\vseq\\
\begin{array}{l}
D_{ij}(\xx^n)+h\GGij(\xx^n)\cdot\VV\geq 0 \hbox{ if } \gammaij^n= 0\vseq\\
D_{ij}(\xx^n)+h\GGij(\xx^n)\cdot\VV = 0 \hbox{ if } \gammaij^n < 0\\
\end{array}

\vseq\\
\hbox{Pairs } (i,k) \hbox{ particle/obstacle :}\vseq\\
\begin{array}{l}
D_{ik}(\xx^n,\yy^{n+1})+h\GGik(\xx^n,\yy^{n+1})\cdot\VV\geq 0 \hbox{ if } \gamma_{ik}^n= 0\vseq\\
D_{ik}(\xx^n,\yy^{n+1})+h\GGik(\xx^n,\yy^{n+1})\cdot\VV = 0 \hbox{ if } \gamma_{ik}^n < 0\\
\end{array}

\end{array}
\right.
\right\}.
$$
\end{rmrk}

\subsection{Finding neighbours}

The most time consuming step in algorithm~\ref{algo_cont_visc_N} is
the projection step~(\ref{algo_N_u}). It is performed using a Uzawa
algorithm which imposes to run matrix/vector products involving the
contacts. However, in order to simulate large collections of particles, it is
essential to avoid loops over the $N(N-1)/2$ possible contacts. To do
so, we notice that it is not necessary to take into account all
contacts at each time-step. Indeed, two particles $i$ and $j$ far enough to each
other at time $t^n$ won't stick at time $t^{n+1}$ and consequently,
the corresponding constraint won't be activated
(ie. $\lambda_{ij}^{n+1}=0$). We denote by $D_{neigh}$ the
distance above which we consider that two particles are not likely to
touch next time-step. Then, the set of pairs of particles one has to
consider at time $t^n$ is:
$$
C_{neigh}(\xx^n)=\left\{ (i,j)\in[1,N]^2\virg i<j \hbox{ and } D_{ij}(\xx^n)\leq D_{neigh}\right\}.
$$

If the pair $(i,j)$ is in the set $C_{neigh}(\xx^n)$, we say that
particles $i$ and $j$ are neighbours. Two particles that are not
neighbours at time $t^n$ won't stick at time $t^{n+1}$ and
consequently, one can restrict the set of constraints at time $t^n$ to:
$$
K_{neigh}(\xx^{n},\ggamma^{n})=\left\{
\begin{array}{l}
\displaystyle\VV\in\IR^{3N} \text{ s.t.}\virg \forall (i,j)\in C_{neigh}(\xx^n)\virg\vseq\\

\hspace{0.5cm}\left|
\begin{array}{l}
\displaystyle D_{ij}(\xx^n)+h\GGij(\xx^n)\cdot\VV\geq 0 \hbox{ if } \gammaij^n\geq 0\vseq\\
\displaystyle D_{ij}(\xx^n)+h\GGij(\xx^n)\cdot\VV = 0 \hbox{ if } \gammaij^n < 0\\
\end{array}
\right.
\end{array}
\right\}.
$$
\begin{rmrk}
This idea not to take into account particles far away from each other is
generally used when considering particles interacting through near
field interaction forces, decreasing with the distance. In that case, it consists in considering
that the force is negligible above a certain distance and consequently, it
is an approximation of the model. In our case, no approximation is
made. Indeed, if $D_{neigh}$ has been chosen sufficiently large, we know
that the pairs of particles that are not belonging to
$C_{neigh}(\xx^n)$ won't interact at time $t^{n+1}$. For example, we
can choose a time
step in order to limit the displacement of the particles to twice
their radius and then set the value of $D_{neigh}$ to a few radiuses.
\end{rmrk}

\begin{figure}[hbtp]
\psfragscanon
\psfrag{l}[r]{\huge $l\nu$}
\psfrag{lp1}[r]{\huge $(l+1)\nu$}
\psfrag{k}[c]{\huge $k\nu$}
\psfrag{kp1}[c]{\huge $(k+1)\nu$}
\psfrag{pi}[l]{\huge particule $i$}
\begin{center}
\resizebox{0.4\textwidth}{!}{\includegraphics{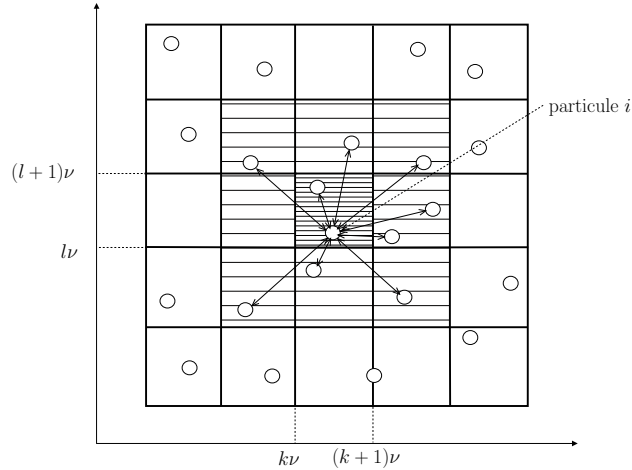}}
\end{center}
\caption{Algorithm to find neighbours: neighbouring boxes and
  distances actually computed.}
\label{fig:voisins}
\end{figure}

To construct $C_{neigh}(\xx^n)$ avoiding the
computation of the $N(N-1)/2$ distances, we choose a bucket
sorting type algorithm. It consists in dividing the computational
domain into boxes of size $\nu>D_{neigh}$ and to compute
distances only for pairs of particles belonging to neighbouring
boxes (see Fig.~\ref{fig:voisins}). 
Note that, because of step~(\ref{algo_N_gamma}), it is not sufficient
to erase at each time-step the former
set of neighbours and to create the new one: one has to transfer the
value of $\gamma_{ij}^n$ if
particles $i$ and $j$ are in contact during these two successive time
steps.

\subsection{Object Oriented Programming Method}

 To build this code,
  we chose to use the object oriented programming method for
  mathematical problems CsiMoon~\cite{CsiMoon}. As a consequence, both
  numerical methods and models can be easily changed. For example, new
  methods  can be chosen and added to the code in order to perform the
 projection step and to construct the set of
  neighbours. This programming method also alows us to take into
  account various models of external environment (dry environment, fluid,
  obstacles of different shapes...), of interparticular interactions
  (cohesion force...) and of contacts (inelastic, gluey model, aggregation...). This leads to a
  modular C++ code SCoPI~\cite{SCoPI}, allowing Simulations of Collections of
  Interacting Particles. This code has already been used to simulate
  gluey particles, crowd motion, wet particles and red-cells (as an
  assembly of rigid particles).

\section{Numerical simulations}

We present in this section numerical simulations of collections of
gluey particles. For visualization reasons, we
only propose here two-dimensional simulations: even though the code is intrinsically
three-dimensional, the motion of the particles is restricted to a vertical
plane. These simulations demonstrate that the algorithm enables to
take great numbers of gluey particles into account. This, together
with section~\ref{sec:couplage}, shows that coupling the gluey
particle algorithm with fluid/particle solvers will make it possible to
simulate dense fluid/particle flows, taking the lubrication force into
account with accuracy.

\subsection{Gluey lotto: influence of roughness}

\begin{figure}[hbtp]
\begin{center}
\resizebox{0.2\textwidth}{!}{\includegraphics{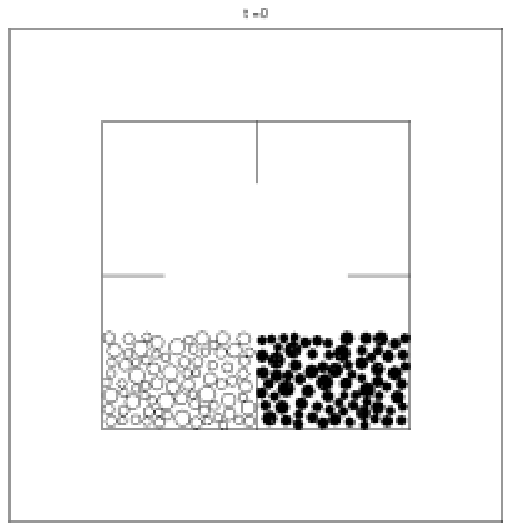}}\hspace{1cm}
\resizebox{0.2\textwidth}{!}{\includegraphics{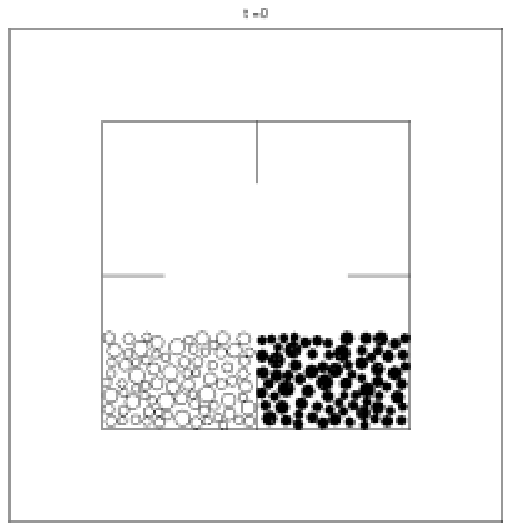}}\hspace{1cm}
\resizebox{0.2\textwidth}{!}{\includegraphics{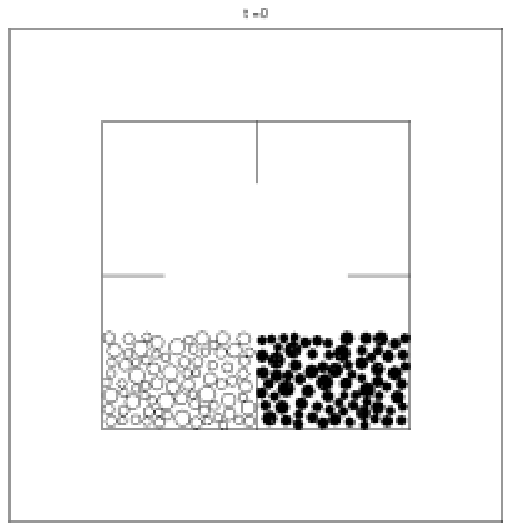}}\\
\resizebox{0.2\textwidth}{!}{\includegraphics{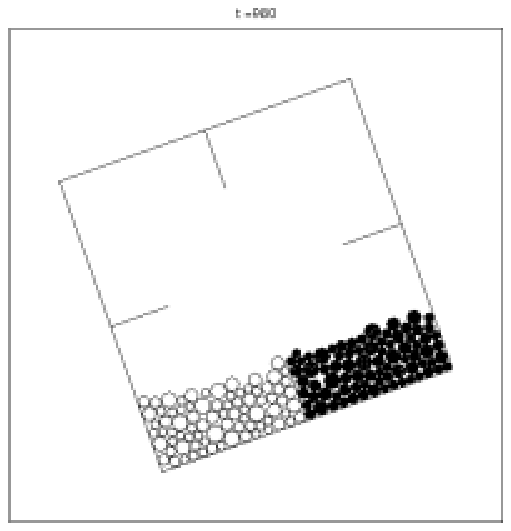}}\hspace{1cm}
\resizebox{0.2\textwidth}{!}{\includegraphics{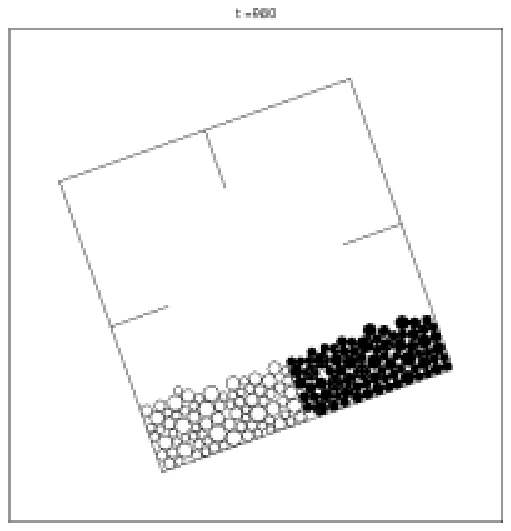}}\hspace{1cm}
\resizebox{0.2\textwidth}{!}{\includegraphics{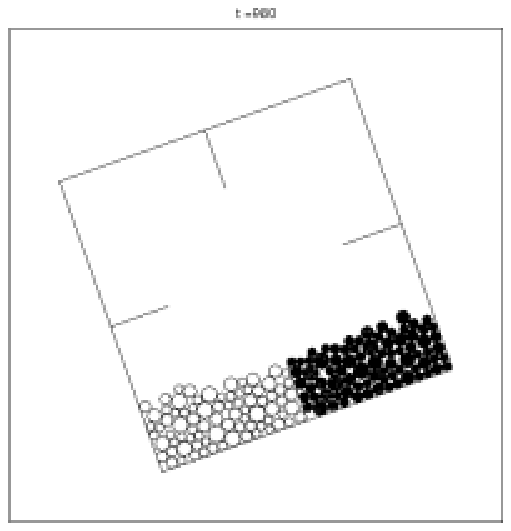}}\\
\resizebox{0.2\textwidth}{!}{\includegraphics{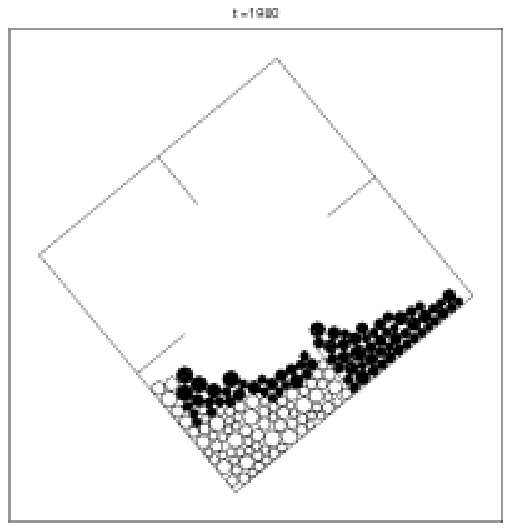}}\hspace{1cm}
\resizebox{0.2\textwidth}{!}{\includegraphics{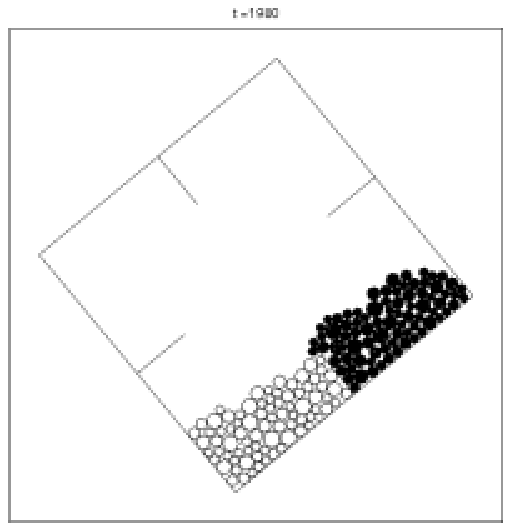}}\hspace{1cm}
\resizebox{0.2\textwidth}{!}{\includegraphics{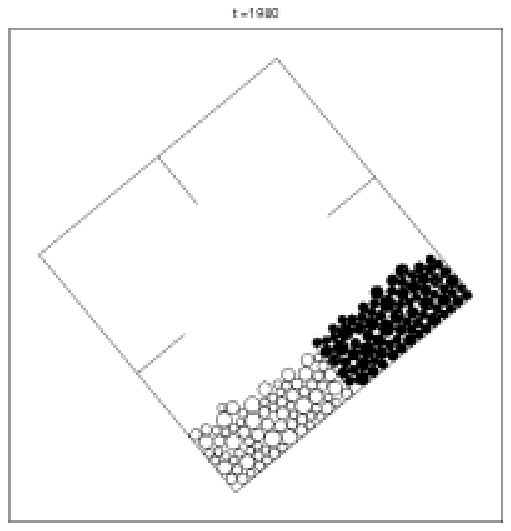}}\\
\resizebox{0.2\textwidth}{!}{\includegraphics{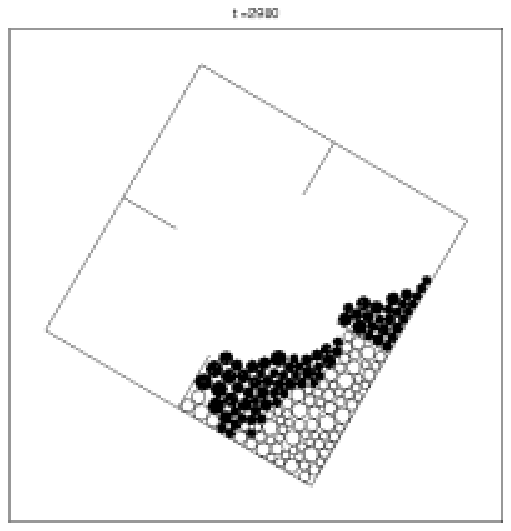}}\hspace{1cm}
\resizebox{0.2\textwidth}{!}{\includegraphics{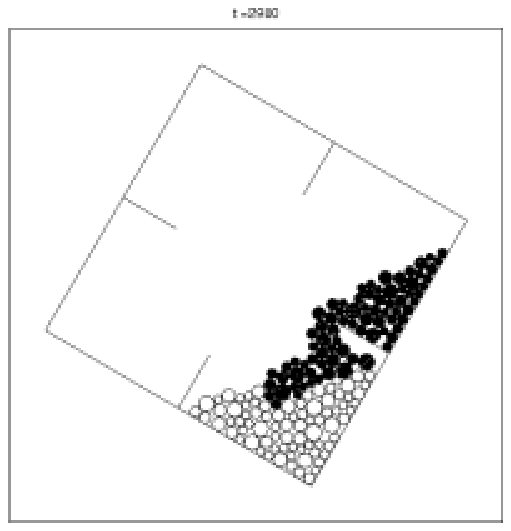}}\hspace{1cm}
\resizebox{0.2\textwidth}{!}{\includegraphics{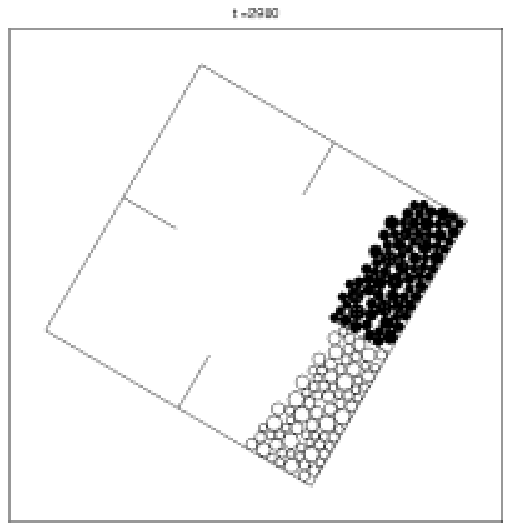}}\\
\resizebox{0.2\textwidth}{!}{\includegraphics{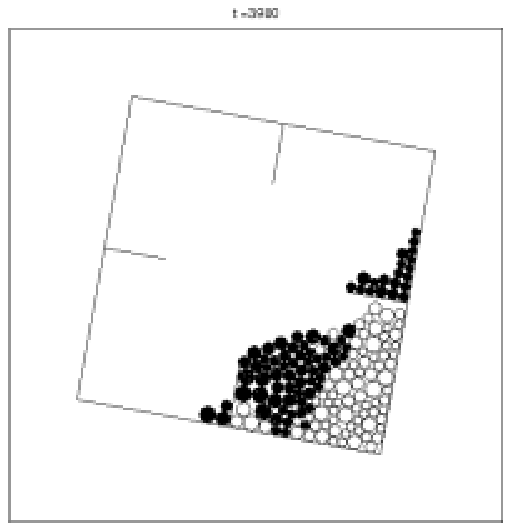}}\hspace{1cm}
\resizebox{0.2\textwidth}{!}{\includegraphics{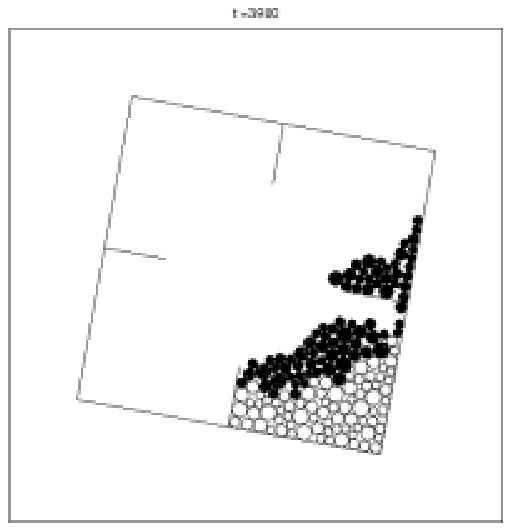}}\hspace{1cm}
\resizebox{0.2\textwidth}{!}{\includegraphics{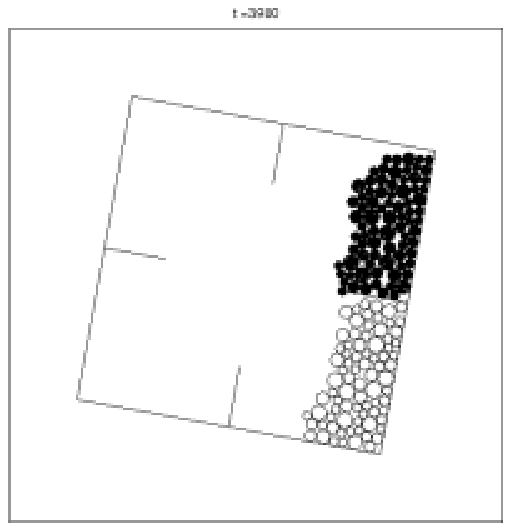}}\\
\resizebox{0.2\textwidth}{!}{\includegraphics{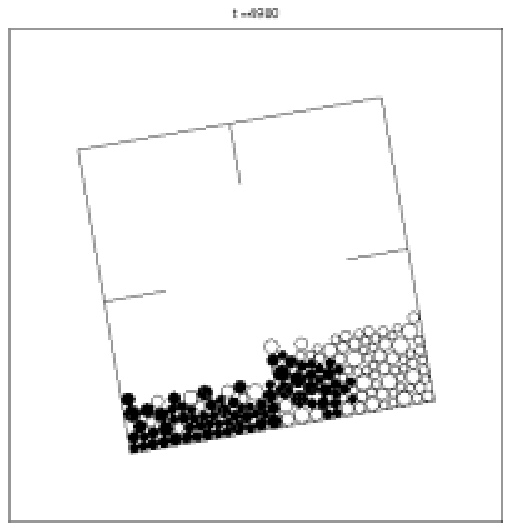}}\hspace{1cm}
\resizebox{0.2\textwidth}{!}{\includegraphics{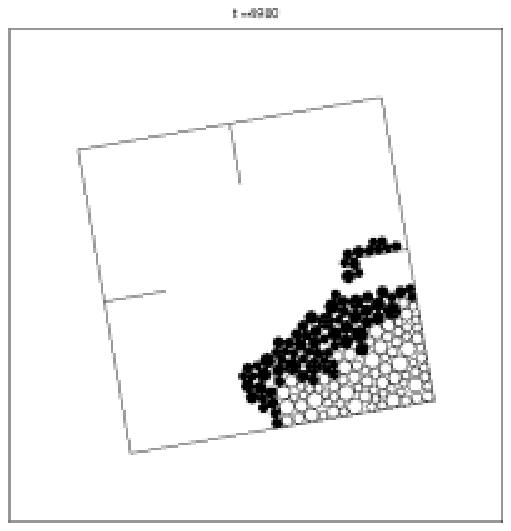}}\hspace{1cm}
\resizebox{0.2\textwidth}{!}{\includegraphics{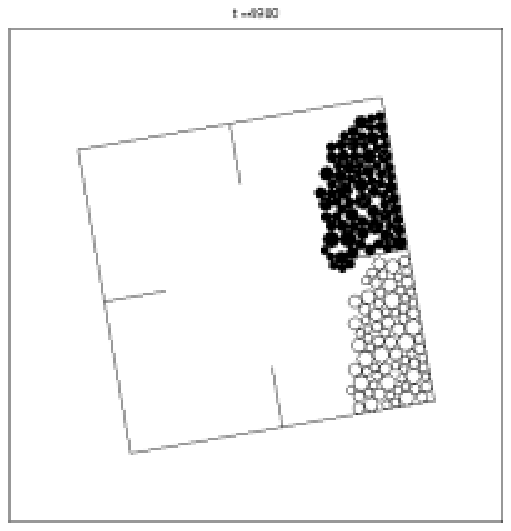}}\\
\end{center}
\caption{Gluey lotto: configurations at
  different time-steps for $\gamma_{min}=0$ on the left (inelastic
contacts), $\gamma_{min}=-1$ in the middle (gluey rough particles) and
$\gamma_{min}=-\infty$ on the right (gluey smooth particles).}
\label{fig_loto_visc}
\end{figure}

The aim of this simulation is to observe the influence of roughness on
the behaviour of multi-particle systems governed by the gluey
particle model. We consider a two-dimensional ``gluey lotto'' made of $160$ particles in
a squared rotating mixer operator. The side lenght of the box is $0.5$
and the radiuses of the particles are taken between $0.007$ and
$0.015$. All particles have the same mass $m=1$ and the gravity
constant $g$ is taken equal to $10$. The $80$ particles initially situated in the left compartment
of the box are black and the $80$ other ones are white. We represent
side by side on figure~\ref{fig_loto_visc} the configurations obtained
at different time steps for  $\gamma_{min}=0$ on the left (inelastic
contacts), $\gamma_{min}=-1$ in the middle (gluey rough particles) and
$\gamma_{min}=-\infty$ on the right (gluey smooth particles). In case
of smooth particles, the  heaps of particles take off from the wall
when they are at the top of the box: as suggested by the
particle/plane model, they take off only when the gravity has balanced
the forces it has itself exerted to push the particles on the bottom
wall. In the rough case, they take off earlier.

%
%

\subsection{Sedimentation of 3000 gluey particles}

\begin{figure}[hbtp]
\begin{center}
\resizebox{0.45\textwidth}{!}{\includegraphics{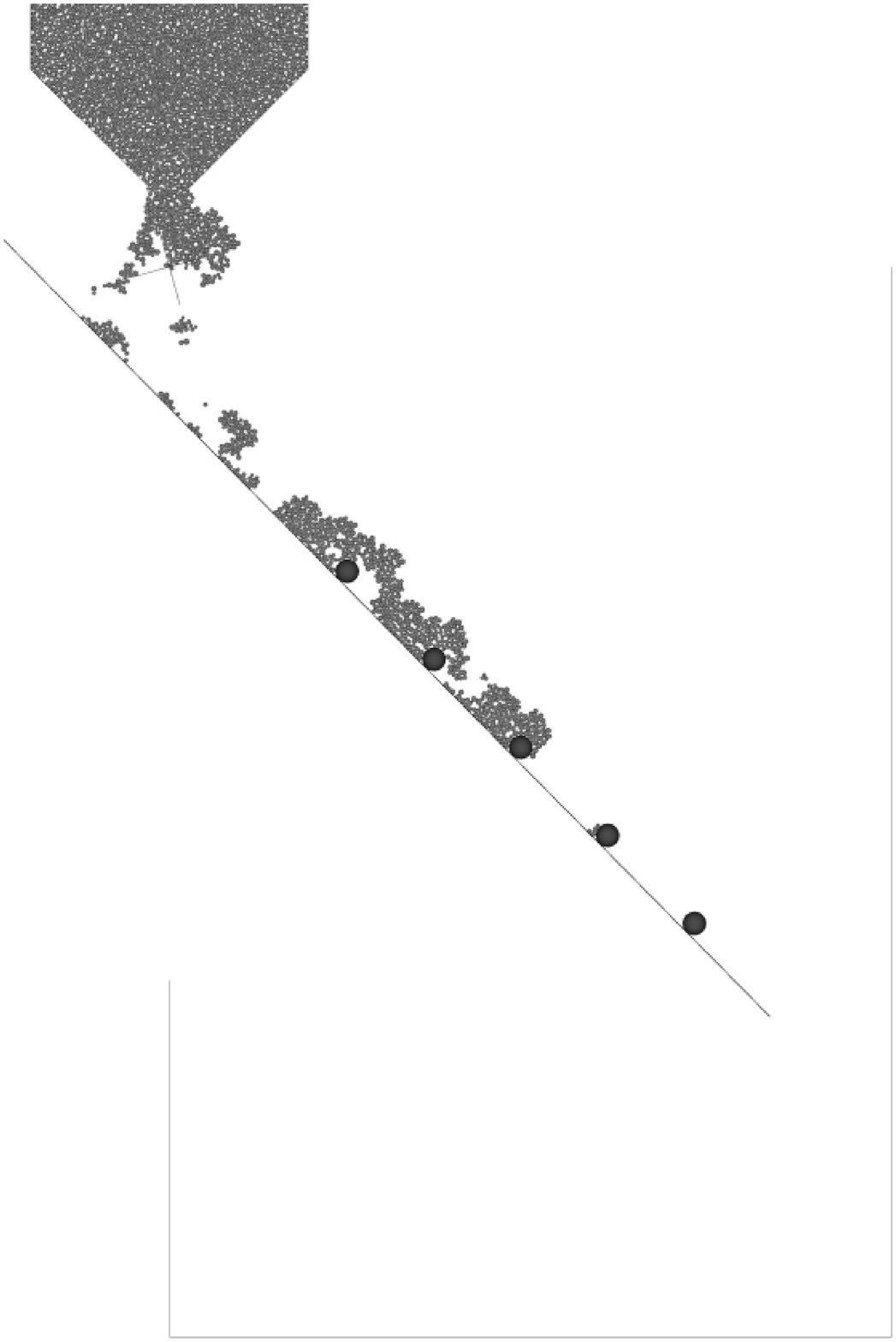}}
\resizebox{0.45\textwidth}{!}{\includegraphics{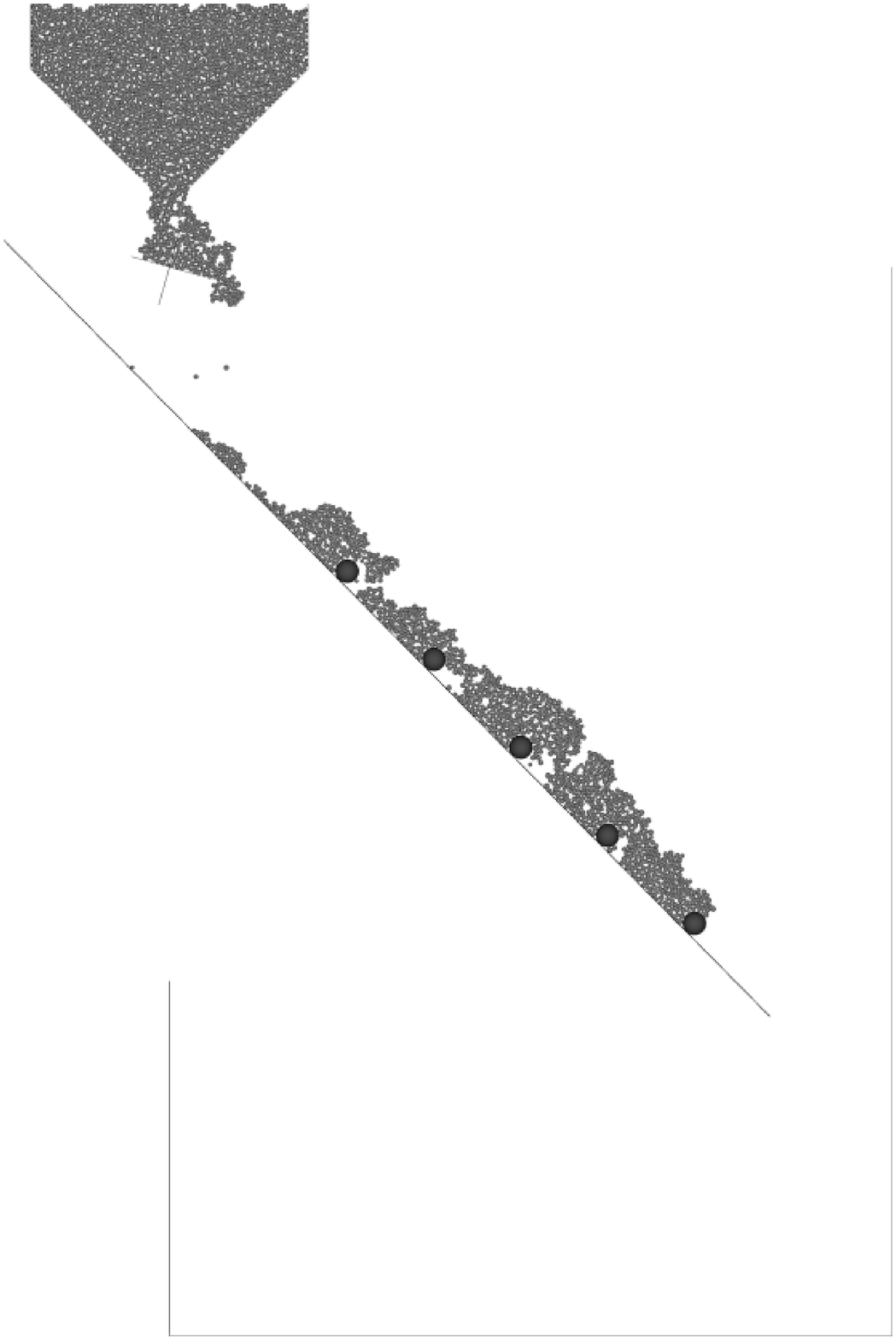}}\\
\resizebox{0.45\textwidth}{!}{\includegraphics{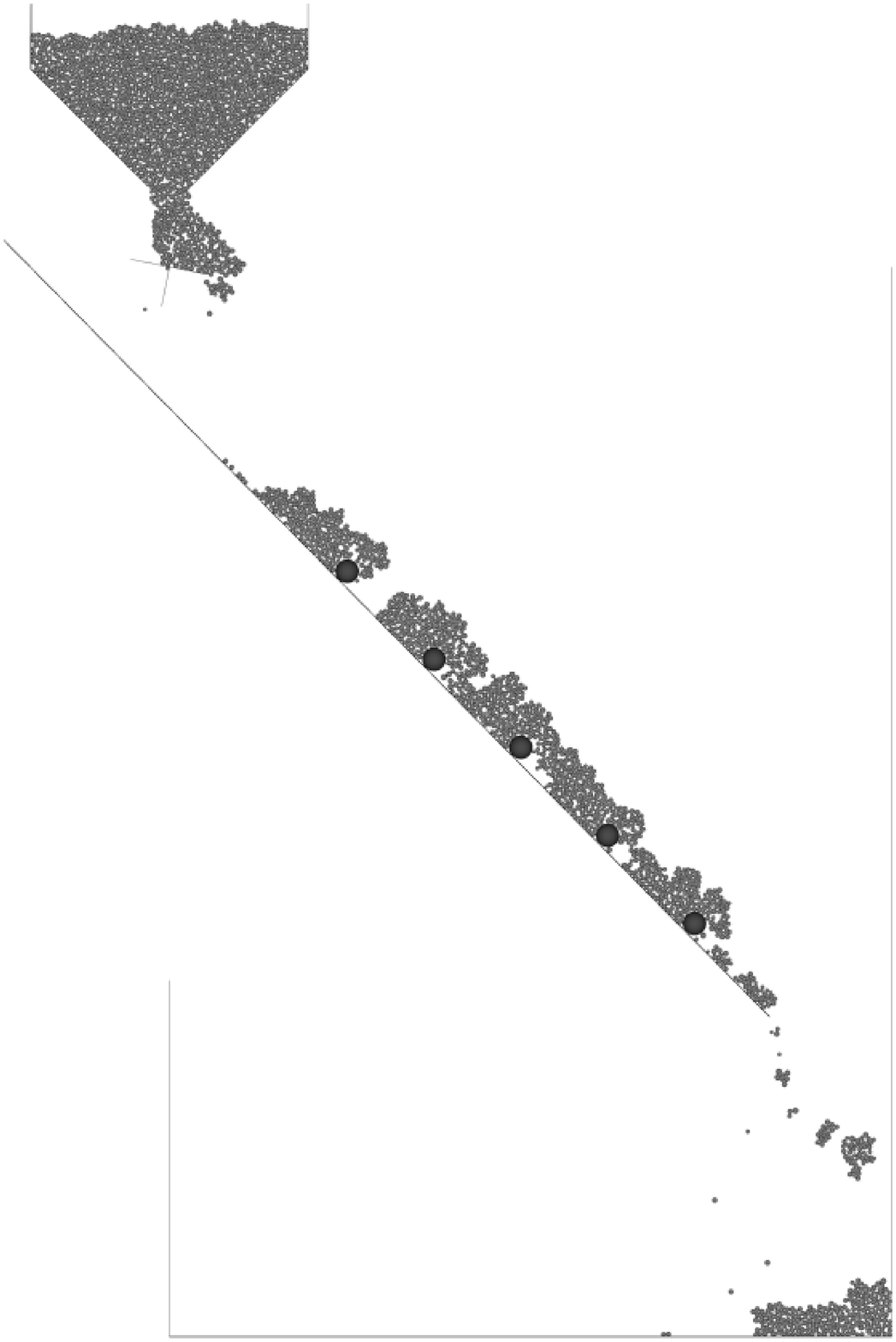}}
\resizebox{0.45\textwidth}{!}{\includegraphics{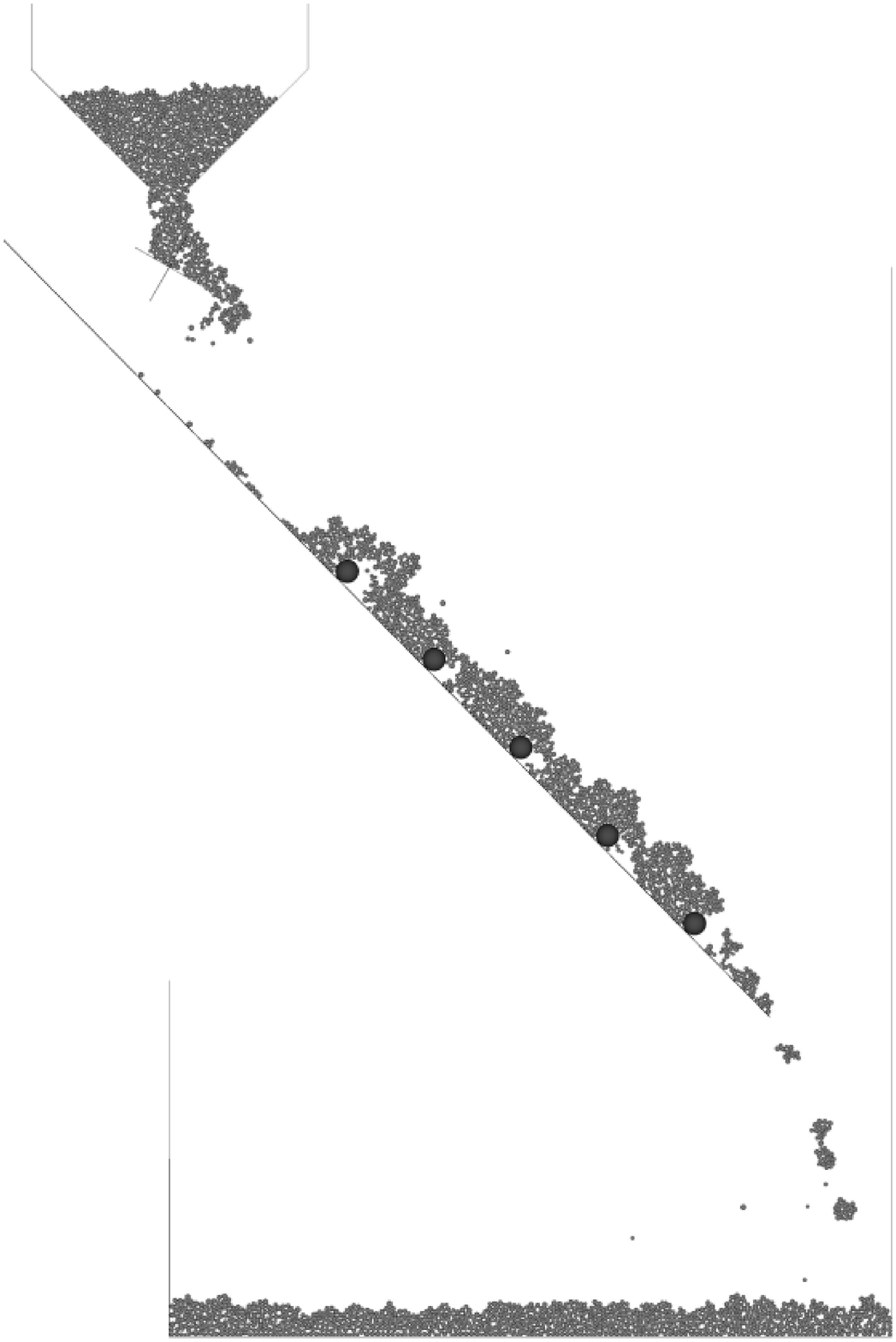}}
\end{center}
\caption{Snapshots of a two-dimensional simulation, 3000 particles:
  configurations at time-steps $n=8826-14226$ (top) and $n=18726-30726$ (bottom).}
\label{fig_SCoPI_2D}
\end{figure}

\begin{figure}[hbtp]
\begin{center}
\resizebox{0.49\textwidth}{!}{\includegraphics{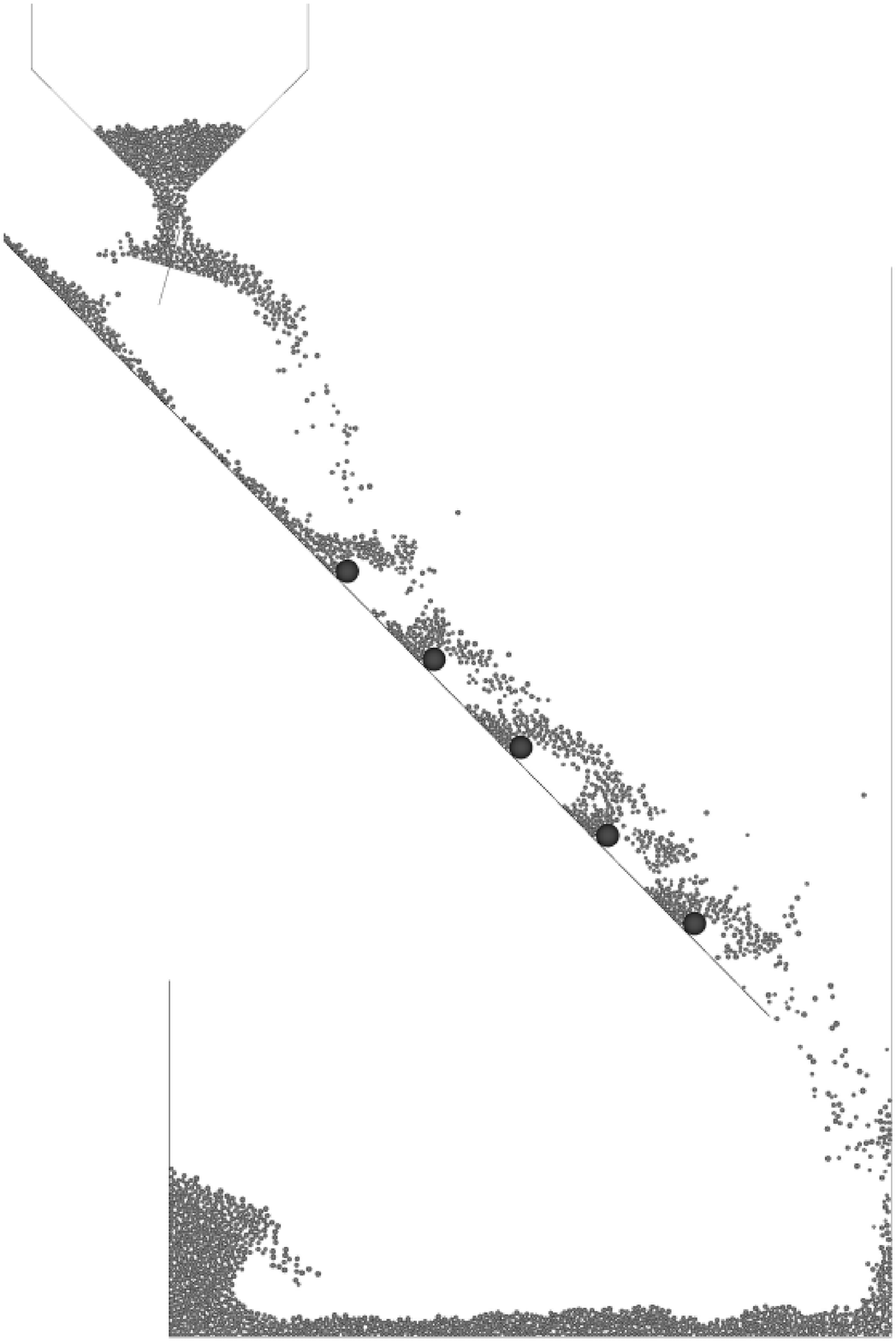}}
\resizebox{0.49\textwidth}{!}{\includegraphics{./IMAGES/IMAGES_SCoPI/SCoPI_2D_bw14226.eps}}
\end{center}
\caption{Snapshot of two-dimensional simulations, 3000 particles: dry
  (left) and viscous (right) simulations at time-step $n=14226$.}
\label{fig_SCoPI_2D_nv}
\end{figure}

\begin{figure}[hbtp]
\begin{center}
\resizebox{0.9\textwidth}{!}{\includegraphics{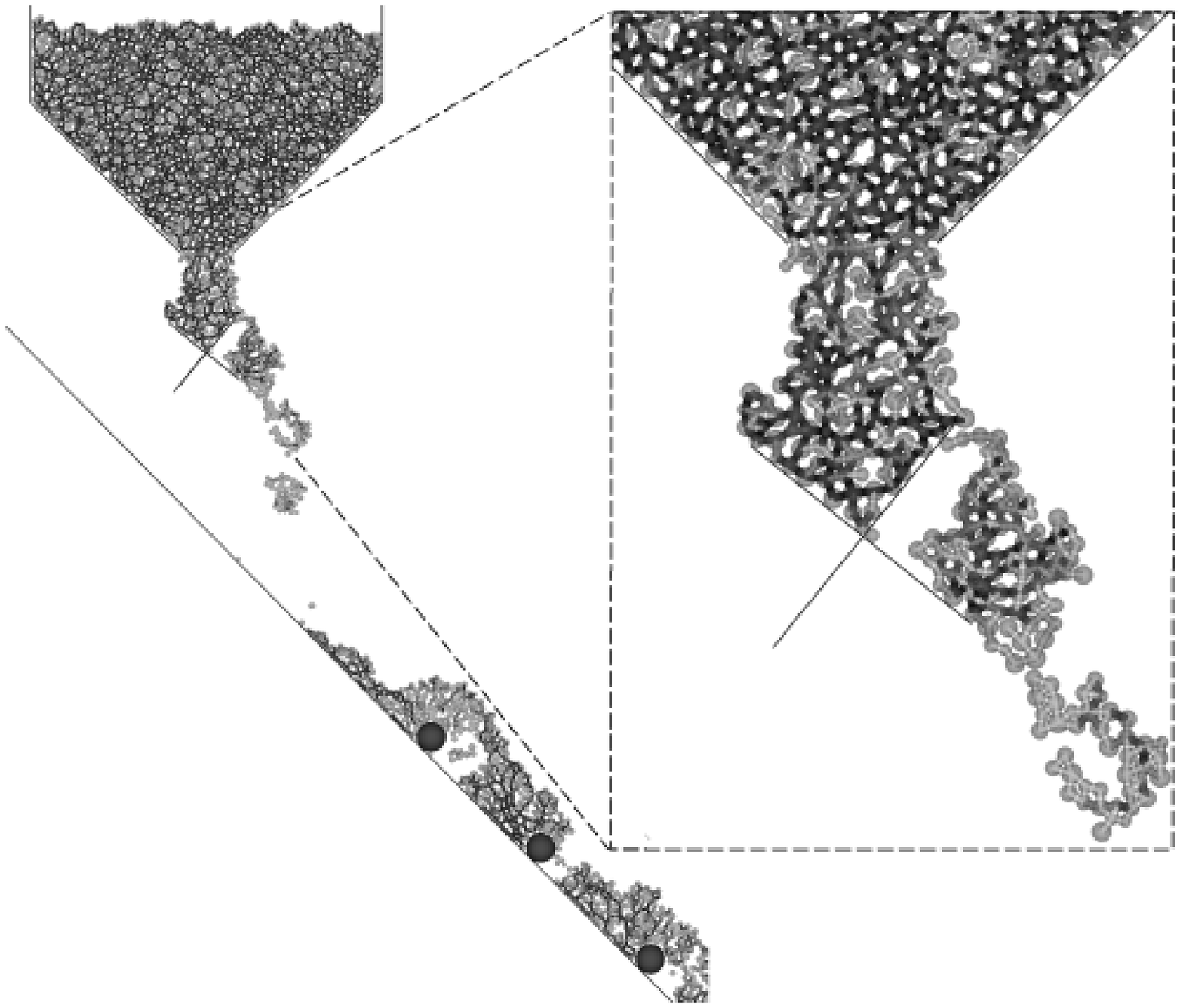}}
\end{center}
\caption{Snapshot of a two-dimensional gluey simulation, 3000 particles:
  configuration and values of $\gamma$ at time-step $n=14226$.}
\label{fig_SCoPI_2D_gamma}
\end{figure}

We consider $3000$ gluey particles sedimenting under gravity with
radiuses between $0.015$ and $0.025$. They are
initially situated above a funnel (random sample of positions) with velocity equal to
zero. All the particles have the same mass $m=2$ and the gravity $g$
is taken equal to $10$. Below the funnel, a wheel rotates around its
axis with angular velocity $\omega=-2$ and throws
the particles on a leaning fixed plane situated below it. Then, the
particles slip along the plane and finally fall in a
container. Some spherical obstacles of radius $r=0.1$ are fixed on the plane to slow the particles
movement. A threshold is imposed on $\gamma$ ($\gamma\geq -10$) to
model roughness. Snapshots of this simulation are presented for different
time-steps on figure~\ref{fig_SCoPI_2D}.
The code also allows us to model dry granular flow involving inelastic contacts. In
figure~\ref{fig_SCoPI_2D_nv} we compare the configurations obtained
at the same time-step for such a simulation and the previous gluey one.
Finally, we plot on figure~\ref{fig_SCoPI_2D_gamma} the values of
$\gamma$ for a given configuration of the gluey simulation. For each contact, a tube is
plotted between the two involved particles and, the larger is
$\gamma_{ij}$ (ie. the more the particles are glued), the more the
grey is dark. We can see the network of the forces leading to a packed configuration in
the funnel. The particles are smoothly unsticking from each other when
leaving the wheel.

\newpage
\text{}

\end{document}